\documentclass[preprint,12pt]{elsarticle}

\usepackage{amssymb}
\usepackage{amsmath}
\usepackage{algorithm}
\usepackage{algorithmic}
\usepackage{bm}
\usepackage{mathrsfs}
\usepackage{makecell}
\usepackage{graphicx}
\usepackage{epsfig}
\usepackage[caption=false,font=normalsize,labelfont=sf,textfont=sf]{subfig}
\usepackage[colorlinks,
linkcolor=red,
anchorcolor=blue,
citecolor=blue
]{hyperref}
\usepackage{tikz}
\usetikzlibrary{spy}

\usepackage{easyReview}

\usepackage{xcolor}
\definecolor{jkcol}{RGB}{45,125,154}
\definecolor{bhcol}{RGB}{247,99,12}
\definecolor{licol}{RGB}{209,52,56}

\newtheorem{thm}{Theorem}

\newtheorem{lem}{Lemma}
\newtheorem{prop}{Proposition}
\newtheorem{defn}{Definition}

\newtheorem{rem}{Remark}

\begin{document}

\begin{frontmatter}


\title{Nonconvex Robust Quaternion Matrix Completion for Imaging Processing\tnoteref{label3}}
\tnotetext[label3]{The work was supported by National Natural Science Foundation of China (12001211, 12071159, 12171168, 12471351) and Natural Science Foundation of Fujian Province, China (2022J01194).}

\author[1]{Baohua Huang}
\affiliation[1]{organization={School of Mathematics and Statistics},
addressline={Fujian Normal University},
city={Fuzhou},
postcode={350117},
country={China}}

\author[2]{Jiakai Chen}
\affiliation[2]{organization={Faculty of Science and Technology},
addressline={United International College (BNU-HKBU)},
city={Zhuhai},
postcode={519087},
country={China}}

\author[3]{Wen Li\fnref{cor1}}
\affiliation[3]{organization={ School of Mathematical Sciences},
addressline={South China Normal University},
city={Guangzhou},
postcode={510631},
country={China}}
\fntext[cor1]{Corresponding author.
E-mail address: liwen@scnu.edu.cn (W. Li).}


%
%
%
%

\begin{abstract}
One of the tasks in color image processing and computer vision is to recover clean data from partial observations corrupted by noise. To this end, robust quaternion matrix completion (QMC) has recently attracted more attention and shown its effectiveness, whose convex relaxation is to minimize the quaternion nuclear norm plus the quaternion $L_1$-norm. However, there is still room to improve due to the convexity of the convex surrogates.
This paper proposes a new nonconvex robust QMC model, in which the nonconvex MCP function and the quaternion $L_p$-norm are used to enhance the low-rankness and sparseness of the low-rank term and sparse term, respectively. An alternating direction method of multipliers (ADMM) algorithm is developed to solve the proposed model and its convergence is given. Moreover, a novel nonlocal-self-similarity-based nonconvex robust quaternion completion method is proposed to handle large-scale data.
Numerical results on color images and videos indicate the advantages of the proposed method over some existing ones.
\end{abstract}

\begin{keyword}
Robust quaternion matrix completion \sep Nonconvex surrogate \sep Nonlocal self-similarity \sep Alternating direction method of multiplier
\end{keyword}

\end{frontmatter}



\section{Introduction}
In the last decade, robust matrix recovery has been widely studied and proven to be very effective in the application of image recovery \cite{jiang2023robust,song2024robust}.
In general, the trick is to stack all the image pixels as column vectors of a matrix, and recovery theories and algorithms are employed to the resulting matrix which is low-rank or approximately low-rank. 
For a color image/video, the traditional matrix-based recovery models are applied to red, green, and blue channels respectively, and may result in color distortion during the recovery process.

As a perfect color image and video representation tool, quaternion has attracted much attention in color image  and video processing \cite{zou2023adaptive, guo2024quaternion}.
By encoding a pixel with RGB channels using a pure quaternion, quaternion-based methods treat a color image/video as a quaternion matrix/tensor.
Compared with matrix/tensor-based methods that need to rearrange the elements, quaternion-based methods can better preserve the color structure of color images/videos.
Therefore, many quaternion-based methods have been proposed and widely used in various applications, such as color image deblurring \cite{huang2022quaternionbased}, color image watermarking \cite{wang2013robust},
color image filtering  \cite{chen2014removing} and color face recognition \cite{zou2016quaternion}.

In the field of color image and video inpainting, Chen et al. \cite{chen2020lowrank} extended the traditional low-rank matrix completion (LRMC) model to the low-rank quaternion matrix completion (LRQMC) model, and employed some nonconvex quaternion rank functions to replace the nuclear norm of a quaternion matrix.
By factorizing a quaternion matrix as two smaller factor quaternion matrices, Miao et al. \cite{miao2020quaternionbased} proposed three LRQMC methods based on quaternion double Frobenius norm, quaternion Frobenius/nuclear norm, and quaternion double nuclear norm, respectively.
Following the work in \cite{miao2020quaternionbased},  Yang et al. \cite{yang2021low} acted the quaternion logarithmic norm (QLN) on two smaller factor quaternion matrices of  the target quaternion matrix, and designed  a new LRQMC algorithm named quaternion logarithmic norm based factorization (QLNF)
for the color image completion.
Based on the low-rank decomposition of the quaternion matrix and quaternion nuclear norm minimization
techniques, Miao et al. \cite{miao2022color} established a LRQMC model for the recovery results.
Yang et al. \cite{yang2021weighted} introduced the quaternion truncated nuclear norm (QTNN) for LRQMC to achieve a more accurate approximation for the quaternion rank function. They also developed the double-weighted quaternion truncated nuclear norm (DWQTNN) to speed up the calculation of QTNN, which adds weighted real diagonal matrices on the residual error quaternion matrix.
In addition to the QMC problem, the quaternion tensor completion (QTC) problem has also received much attention. With the help of Tucker rank, the global low-rank prior
to quaternion tensor is encoded as the nuclear norm of unfolding quaternion matrices, and then a QTC algorithm was derived for the color image  and video recovery  \cite{miao2020lowrank}.
Later, Miao et al. \cite{miao2024quaternion} proposed a low rank quaternion tensor completion (LRQTC)   model for the color image and video inpainting which adopts the quaternion weighted nuclear norm (QWNN) of  mode-$n$ canonical unfolding quaternion matrices to capture the global low quaternion tensor train (QTT) rank, and the $L_1$-norm of a quaternion tensor in a transformed domain to capture the sparseness.

The above mentioned LRQMC model mainly aims at the color images  under a low sampling ratio without noise corruption. For reconstructing low-rank matrices from incomplete and corrupted observations,
Jia et al. \cite{jia2019robust} proposed a robust quaternion matrix completion (robust QMC) model. Concretely, the robust QMC is to minimize a hybrid optimization problem involving both the quaternion nuclear norm (QNN) of the low-rank part and the quaternion $L_1$-norm of the sparse part under the limited sample constraints, i.e.,
\begin{equation}\label{1.1}
\min\limits_{\mathbf{L},\mathbf{S}}\|\mathbf{L}\|_*+
\lambda \|\mathbf{S}\|_1,~~\text{s.t.}~~
\mathcal{P}_\Omega(\mathbf{L}+\mathbf{S})=  \mathcal{P}_\Omega(\mathbf{X}),
\end{equation}
where $\mathbf{X}\in\mathbb{Q}^{n_1\times n_2}$ is a noisy quaternion matrix, $\mathbf{L}\in\mathbb{Q}^{n_1\times n_2}$ is the target low-rank quaternion matrix,
$\mathbf{S}\in\mathbb{Q}^{n_1\times n_2}$ is a sparse quaternion matrix and acts as the corrupted  data, $\lambda$ is a regularization parameter used to balance the low rank and sparse parts, $\Omega$ is an index set and $\mathcal{P}_\Omega$ is the unitary projection onto $\Omega$ such that the entries in the set $\Omega$ are given while the remaining entries are missing. 
The robust QMC model (\ref{1.1}) performs quite well in the color image completion problem. However, the gap between the rank function and QNN may be large, especially when some of the singular values of the original quaternion matrix are very large. Besides, the quaternion $L_1$-norm  is a coarse estimation of the sparse part and leads some of the shortcomings in theory and experiment.

In order to better approximate the real matrix rank function, a lot of nonconvex real matrix rank   surrogates have
been proposed, such as,  the schatten $q$-norm \cite{frank1993statistical},  weighted schatten $q$-norm \cite{xie2016weighted}, weighted nuclear norm \cite{gu2017weighted}, and log-determinant penalty  \cite{kang2015logdet}. Along this way, many nonconvex quaternion matrix rank surrogates were proposed by Chen et al. \cite{chen2020lowrank}  in the quaternion matrix approximate problem,   e.g., Laplace, Geman, and weighted schatten $q$. Unlike the over-penalty of nuclear norm for large singular values,  nonconvex penalties of singular values can result in better approximation of the quaternion matrix rank \cite{chen2020lowrank}. This motivates us to  seek a
more accurate approximation of quaternion matrix rank under the nonconvex setting.
It is known that QNN is the $L_1$-norm of all singular values of a quaternion matrix. Now that the convex regularization of   QNN  achieves more accurate approximation of quaternion rank, we have reason to believe that the similar case also arises for the quaternion $L_1$-norm. It will be interesting to study the nonconvex penalty for sparse estimation instead of using the convex $L_1$-penalty.

Using the low-rankness of underlying quaternion data and the sparseness of sparse corruptions, this paper proposes and develops a nonconvex robust low-rank quaternion matrix completion model with nonconvex regularization, which aims to recover a quaternion matrix corrupted by sparse noise with partial observations. More concretely, we first choose a nonconvex quaternion matrix rank approximation related to the minimax concave penalty (MCP) function
which is continuous and unbiased \cite{fan2001variable}.
Second, we employ the penalty function $f(x)=|x|^p (0<p<1)$ on every entry of sparse corruptions. The nonconvex functions are performed on the singular values of the low-rank quaternion matrix and all entries of sparse corruptions, respectively, which are beneficial for promoting the low-rankness of underlying quaternion matrix and enhancing the sparsity of sparse corruptions better compared with QNN and quaternion $L_1$-norm.

It is known that there are a lot of repeated local patterns across a natural image, which is called nonlocal self-similarity (NSS). The nonlocal strategy for image processing was already discussed in \cite{dabov2007image}. Recently, Jia et al. \cite{jia2022nonlocal} proposed a NSS-based QMC (TNSS-based QMC) to recover color images/videos from incomplete and corrupted entries. The TNSS-based QMC searches similar patches of a color image/video to generate a low-rank quaternion sub-matrix and applies QMC for each sub-matrix. Compared with the original data, each resulting sub-matrix would have a lower rank. This paper develops an NSS-based robust low-rank quaternion matrix completion model using nonconvex regularization to better approximate the matrix rank.

The contributions of this paper are summarized as follows.
\begin{itemize}
\item We propose a nonconvex approach for robust quaternion matrix completion (NRQMC). In particular, the MCP function is used as a nonconvex surrogate of matrix rank for a more accurate approximation, and the quaternion $L_p$-norm is developed for better capturing the sparseness of the noisy component.
\item We establish an alternating direction method of multipliers (ADMM) algorithm to solve the proposed NRQMC model. The thresholding operator for the MCP penalty is used to solve the low-rank part subproblem. The quaternion $L_p$ thresholding technique is used to solve the quaternion $L_p$-norm minimization problem.
\item The NSS prior based on quaternion representation is applied and a new NSS-based nonconvex robust quaternion completion method (NRQMC-NSS) is developed for the large-scale color image and video recovery tasks.
\item Numerical experiments on color images and videos demonstrate that the proposed model and algorithm can recover the color and geometric properties such as color edges and textures, and the numerical performances are better than some existing methods.
\end{itemize}

This paper is organized as follows.
Section \ref{sec:pre} reviews some useful notations and definitions.
In Section \ref{sec3}, we develop a nonconvex low-rank approximation based on the MCP function and the $L_p$-norm for quaternion matrices and propose a novel nonconvex robust quaternion matrix completion (NRQMC) model. An ADMM algorithm is established to solve the proposed model and its convergence is discussed.
In Section \ref{sec4}, we introduce the proposed NRQMC-NSS model for large-scale color image and video recovery tasks.
In Section \ref{sec5}, numerical results are presented to show the superiority of the proposed models compared with some existing models. This paper ends with some concluding remarks in Section \ref{sec:con}.

\section{Preliminaries}
\label{sec:pre}
In this section, we introduce some definitions and notations and review some results about quaternion matrix optimization. Let $\mathbb{R}$ and $\mathbb{C}$ be the sets of real numbers and complex numbers, respectively.
Quaternion, invented by Hamilton \cite{hamilton2010elements},  has a real part and three imaginary parts given by
\begin{equation*}
\mathbf{q}=q_0+q_1\mathbf{i}+q_2\mathbf{j}+q_3\mathbf{k},
\end{equation*}
where $q_0, q_1, q_2, q_3\in\mathbb{R}$,  and  three imaginary units $\mathbf{i}$, $\mathbf{j}$, $\mathbf{k}$ obey the following quaternion rules
\begin{equation*}
\mathbf{i}^2=\mathbf{j}^2=\mathbf{k}^2=-1,
~\mathbf{i}\mathbf{j}=-\mathbf{j}\mathbf{i}=\mathbf{k},
~\mathbf{j}\mathbf{k}=-\mathbf{k}\mathbf{j}=\mathbf{i},
~\mathbf{k}\mathbf{i}=-\mathbf{i}\mathbf{k}=\mathbf{j}.
\end{equation*}
Here a boldface symbol indicates that it is a quaternion number, vector, or matrix. Denote $\mathbb{Q}$ the space with quaternion numbers. The quaternion $\mathbf{q}\in \mathbb{Q}$ can be rewritten as
\begin{equation*}
\mathbf{q}=\mathrm{S} \mathbf{q} +\mathrm{V} \mathbf{q},
\end{equation*}
where the scalar (real) part  is denoted by $q_0=\mathrm{S} \mathbf{q}=\mathfrak{R}(\mathbf{q})$, whereas the vector (imaginary) part $\mathrm{V} \mathbf{q}=\mathfrak{I}(\mathbf{q})=q_1\mathbf{i}+q_2\mathbf{j}+q_3\mathbf{k}$ comprises the three imaginary parts. If $\mathfrak{R}(\mathbf{q})=0$, then $\mathbf{q}$ is called a pure quaternion.
The conjugate of a quaternion $\mathbf{q}$  is defined as $\mathbf{q}^*=q_0-q_1\mathbf{i}-q_2\mathbf{j}-q_3\mathbf{k}$,
while the conjugate of the product satisfies $(\mathbf{p}\mathbf{q})^*=\mathbf{q}^*\mathbf{p}^*$.
The modulus of a quaternion is defined
as $|\mathbf{q}|=\sqrt{\mathbf{q}\mathbf{q}^*}=\sqrt{q_0^2+q_1^2+q_2^2+q_3^2}$ and it holds
$|\mathbf{p}\mathbf{q}|=|\mathbf{p}||\mathbf{q}|$.
An $n_1\times n_2$ quaternion matrix is of the form
\begin{equation*}
\mathbf{X}=X_0+X_1\mathbf{i}+X_2\mathbf{j}+X_3\mathbf{k},
\end{equation*}
where $X_0,X_1,X_2,X_3\in\mathbb{R}^{n_1\times n_2}$.
Let $\mathbb{Q}^{n_1\times n_2}$ denote the set of all $n_1\times n_2$ quaternion matrices.
A pure quaternion matrix is a matrix whose elements are pure
quaternions ($X_0 = 0$) or zero.

In the RGB color space, every
pixel can be expressed as a pure quaternion, $r\mathbf{i}+g\mathbf{j}+b\mathbf{k}$,
where $r, g, b$ are the values of red, green, and blue
components, respectively.
An $n_1\times n_2$  color image    can be represented by an $n_1\times n_2$ quaternion matrix $\mathbf{X}$ with the following form
\begin{equation}\label{2.1}
\mathbf{X}_{ij}=R_{ij}\mathbf{i}+G_{ij}\mathbf{j}+B_{ij}\mathbf{k},~~1\leq i\leq n_1, 1\leq j\leq n_2,
\end{equation}
where $R_{ij}$, $G_{ij}$, and $B_{ij}$ are the red, green, and blue pixel values, respectively, at the location $(i,j)$ in the image.

For     $\mathbf{X}=(\mathbf{x}_{ij}) \in\mathbb{Q}^{n_1\times n_2}$, let $\mathbf{X}^*=(\mathbf{x}_{ji}^*)\in\mathbb{Q}^{n_2\times n_1}$ be    the conjugate transpose of $\mathbf{X}$.  The unit quaternion matrix $\mathbf{I}$ is just as the classical unit matrix. For a square quaternion matrix $\mathbf{X}\in\mathbb{Q}^{n\times n}$, we say $\mathbf{X}$ is   unitary if $\mathbf{X}^*\mathbf{X} = \mathbf{X}\mathbf{X}^* = \mathbf{I}$.     The sum of all the diagonal elements of $\mathbf{X}\in\mathbb{Q}^{n\times n}$ is called the trace of $\mathbf{X}$, denoted by $\mathrm{tr}(\mathbf{X})$.
We use   $\mathrm{rank}(\mathbf{X})$ to denote the
rank of $\mathbf{X}\in \mathbb{Q}^{n_1\times n_2}$ which is  the maximum number of right   linearly independent columns of a quaternion matrix $\mathbf{X}$.
For $\mathbf{X}=(\mathbf{x}_{ij})\in\mathbb{Q}^{n_1\times n_2}$, $\mathrm{absQ}(\mathbf{X})=(|\mathbf{x}_{ij}|)\in\mathbb{Q}^{n_1\times n_2}$ and $\mathrm{signQ}(\mathbf{x}_{ij})=\mathbf{x}_{ij}/|\mathbf{x}_{ij}|$ (if $|\mathbf{x}_{ij}|\neq 0$) or 0 (otherwise).

The singular value decomposition (QSVD) of a quaternion matrix $\mathbf{X}\in \mathbb{Q}^{n_1\times n_2}$ is given by
\begin{equation}\label{2.2}
\mathbf{X}=\mathbf{U}\mathbf{\Sigma}\mathbf{V}^*,
\end{equation}
where $\mathbf{\Sigma}=\mathrm{diag}(\sigma_1,\cdots,\sigma_{\min\{n_1,n_2\}})\in\mathbb{R}^{n_1\times n_2}$ with $\sigma_i\geq 0$, and $\mathbf{U}\in\mathbb{Q}^{n_1\times n_1}$ and $\mathbf{V}\in\mathbb{Q}^{n_2\times n_2}$ are two unitary quaternion matrices \cite{zhang1997quaternions}. For convenience, we always assume that $n_{(1)}=\max\{n_1,n_2\}$ and $n_{(2)}=\min\{n_1,n_2\}$.

The following quaternion vector and matrix norms are used in the sequel.

\begin{defn}\label{defn2.1}
{\rm (1)} Let $\mathbf{a}=(\mathbf{a}_i)\in\mathbb{Q}^n$ be a quaternion vector. Then the $L_1$-norm $\|\mathbf{a}\|_1=\sum\limits_{i=1}^n|\mathbf{a}_i|$, the $2$-norm $\|\mathbf{a}\|_2=\sqrt{\sum\limits_{i=1}^n|\mathbf{a}_i|^2}$ and the infinity norm $\|\mathbf{a}\|_\infty=\max_{1\leq i \leq n}|\mathbf{a}_i|$.

{\rm (2)} Let $\mathbf{A}=(\mathbf{a}_{ij})\in\mathbb{Q}^{n_1\times n_2}$ be a quaternion matrix. Then the $L_1$-norm $\|\mathbf{A}\|_1=\sum\limits_{i=1}^{n_1}\sum\limits_{j=1}^{n_2}|\mathbf{a}_{ij}|$,  the Frobenius norm $\|\mathbf{A}\|_F=\sqrt{\sum\limits_{i=1}^{n_1}\sum\limits_{j=1}^{n_2}|\mathbf{a}_{ij}|^2}
$,  the infinity norm $\|\mathbf{A}\|_\infty=\max_{i,j}|\mathbf{a}_{ij}|$, the spectral norm $\|\mathbf{A}\|_2=\max\{\sigma_1,\cdots,\sigma_r\}$ and the nuclear norm $\|\mathbf{A}\|_*=\sum\limits_{i=1}^{r}\sigma_i$, where $\sigma_1,\cdots,\sigma_r$ are nonzero singular values of $\mathbf{A}$.
\end{defn}

In the vector space $\mathbb{Q}^{n_1\times n_2}$, the real inner product is given by \cite{qi2020quaternion}
\begin{equation}\label{new2.1}
\langle \mathbf{X}, \mathbf{Y}\rangle = \mathrm{Re}(\mathrm{tr}(\mathbf{X}^*\mathbf{Y}))
\end{equation}
for all $\mathbf{X}, \mathbf{Y}\in \mathbb{Q}^{n_1\times n_2}$.
The norm of a quaternion matrix generated by this inner product space is denoted by $\|\cdot\|$. Then,
for $\mathbf{X}\in \mathbb{Q}^{n_1\times n_2}$, we have
\begin{equation*}
\|\mathbf{X}\|^2 = \langle \mathbf{X}, \mathbf{X}\rangle = \mathrm{Re}(\mathrm{tr}(\mathbf{X}^* \mathbf{X})) = \mathrm{tr}(\mathbf{X}^* \mathbf{X}) = \|\mathbf{X}\|^2_F
\end{equation*}
since $\mathrm{tr}(\mathbf{X}^* \mathbf{X})$ is real.
Direct  calculations also give
\begin{equation}\label{new2.2}
\langle \mathbf{X}, \mathbf{Y}\mathbf{Z}\rangle = \langle \mathbf{Y}^*\mathbf{X},\mathbf{Z}\rangle = \langle \mathbf{X}\mathbf{Z}^*, \mathbf{Y}\rangle
\end{equation}
and
\begin{equation}\label{neq2.3}
\|\mathbf{Q}\mathbf{E}\|_F=\|\mathbf{E}\mathbf{Q}\|_F=\|\mathbf{E}\|_F,
\end{equation}
where $\mathbf{X}\in \mathbb{Q}^{n_1\times n_2}$, $\mathbf{Y}\in \mathbb{Q}^{n_1\times l}$,  $\mathbf{Z}\in \mathbb{Q}^{l\times n_2}$,
$\mathbf{E}\in \mathbb{Q}^{n_2\times n_2}$,  and $\mathbf{Q}\in \mathbb{Q}^{n_2\times n_2}$ is an unitrary quaternion matrix.

Let $\mathscr{S}$ be a finite dimensional  space, we say a function $f:\mathscr{S}\rightarrow [-\infty,+\infty]$ is proper if $f(x)<+\infty$ for at least one $x\in\mathscr{S}$, and $f(x)>-\infty$ for all $x\in\mathscr{S}$ \cite{mordukhovich2020variational}. For a proper and lower semicontinuous function $f:\mathscr{S}\rightarrow (-\infty,+\infty]$, the proximal mapping associated with $f$ at $y$ is specified by
\begin{equation*}
\mathrm{Prox}_f(y)=\arg\min\limits_{x\in\mathscr{S}}\Big\{f(x)+\dfrac{1}{2}\|x-y\|^2\Big\},~\forall
y\in \mathscr{S}.
\end{equation*}
For a nonconvex function $f:\mathbb{R}^n\rightarrow (-\infty,+\infty]$, the subdifferential of $f$ at $x$ \cite[Definition 8.3]{mordukhovich2020variational}, denoted as $\partial f(x)$, is given by
\begin{equation*}
\partial f(x)=\{y\in\mathbb{R}^n:\exists x_k\rightarrow x,f(x_k)\rightarrow f(x),y_k\rightarrow y~\text{with}~y_k\in f'(x_k)~\text{as}~k\rightarrow+\infty\},
\end{equation*}
where $f'(x_k)$ denotes the Fr\'{e}chet differential of $f$ at $x$.

\begin{defn}[\cite{qi2020quaternion,chen2022color}]\label{defn2.2}
Let $h:\mathbb{Q}^{n_1\times n_2} \rightarrow\mathbb{R}$. We say $h$ is differentiable at $\mathbf{X}=X_0+X_1\mathbf{i}+X_2\mathbf{j}+X_3\mathbf{k}$ if $\dfrac{\partial h}{\partial X_i}$ exists at $X_i$ for $i=0,1,2,3$, and the gradient of $h$ is specified by
\begin{equation}\label{2.3}
	\nabla h(\mathbf{X})=\dfrac{\partial h}{\partial X_0}+
	\dfrac{\partial h}{\partial X_1}\mathbf{i}+
	\dfrac{\partial h}{\partial X_2}\mathbf{j}+
	\dfrac{\partial h}{\partial X_3}\mathbf{k}.
\end{equation}
If $\dfrac{\partial h}{\partial X_i}$ exists in a neighborhood of $X_i$, and it is continuous at $X_i$ for $i=0,1,2,3$, then we say $h$ is continuous differentiable at $\mathbf{X}$. If $h$ is continuous differentiable for any $\mathbf{X}\in\mathbb{Q}^{n_1\times n_2}$, then we sat $g$ is continuously differentiable.
\end{defn}

The directional derivative of $h$ at $\mathbf{X}\in\mathbb{Q}^{n_1\times n_2}$ in the direction $\Delta\mathbf{X}\in\mathbb{Q}^{n_1\times n_2}$ is defined by
\begin{equation*}
h'(\mathbf{X};\bm{\Delta}\mathbf{X})=
\lim\limits_{t\rightarrow 0,t\in\mathbb{R}}\dfrac{h(\mathbf{X}+t\bm{\Delta}\mathbf{X})-h(\mathbf{X})}{t}.
\end{equation*}
We see that the directional derivative of $h$ is real, while the gradient of $h$ is in $\mathbb{Q}^{n_1\times n_2}$. However, they can be connected via the real inner product (\ref{new2.1}).

\begin{prop}[\cite{qi2020quaternion}]\label{prop2.1}
For $\mathbf{X},\bm{\Delta}\mathbf{X}\in\mathbb{Q}^{n_1\times n_2}$, it holds
$h'(\mathbf{X};\bm{\Delta}\mathbf{X})= \langle\nabla h(\mathbf{X})$, $\bm{\Delta}\mathbf{X}\rangle$.
If $h$ is continuous differentiable, then
$h(\mathbf{X}+\bm{\Delta}\mathbf{X})=h(\mathbf{X})+\langle\nabla h(\mathbf{X})$, $\bm{\Delta}\mathbf{X}\rangle
+o(\|\bm{\Delta}\mathbf{X}\|_F)$.
\end{prop}

For a quaternion matrix $\mathbf{X} = X_0+X_1\mathbf{i}+X_2\mathbf{j}+X_3\mathbf{k}\in\mathbb{Q}^{n_1\times n_2}$, let $\mathscr{R}(\mathbf{X})=(X_0,X_1,X_2,X_3)\in\mathbb{R}^{n_1\times n_2}\times \mathbb{R}^{n_1\times n_2}\times \mathbb{R}^{n_1\times n_2}\times \mathbb{R}^{n_1\times n_2}$. Let $h:\mathbb{Q}^{n_1\times n_2} \rightarrow\mathbb{R}$. Then $h$ can be regarded as a function of $\mathscr{R}(\mathbf{X})$ and denote such function as $h^R$. We say $h:\mathbb{Q}^{n_1\times n_2} \rightarrow\mathbb{R}$ is twice continuously differentiable  with respect to $\mathbf{X}$ if $h^R$ is twice continuously differentiable with respect to $\mathscr{R}(\mathbf{X})$ \cite{qi2020quaternion}. When $h$ is twice continuously differentiable, we consider the second order derivative  $\nabla^2 h(\mathbf{X})$. As in \cite{qi2020quaternion}, it is more convenient to consider  $\nabla^2 h(\mathbf{X})\bm{\Delta}\mathbf{X}$ and $\langle\nabla^2 h(\mathbf{X})\bm{\Delta}\mathbf{X},\bm{\Delta}\mathbf{X}\rangle$.
Suppose that
\begin{equation*}
\nabla h(\mathbf{X}+\bm{\Delta}\mathbf{X})-\nabla h(\mathbf{X})=\vartheta(\mathbf{X}+\bm{\Delta}\mathbf{X})
+\varsigma(\mathbf{X}+\bm{\Delta}\mathbf{X}),
\end{equation*}
where $\vartheta$ is real linear in $\bm{\Delta}\mathbf{X}$ in the sense that for any $a,b\in\mathbb{R}$ and
$\bm{\Delta}\mathbf{X}^{(1)}$, $\bm{\Delta}\mathbf{X}^{(2)}\in\mathbb{Q}^{n_1\times n_2}$,
\begin{equation*}
\vartheta(\mathbf{X}+a\bm{\Delta}\mathbf{X}^{(1)}+b\bm{\Delta}\mathbf{X}^{(1)})
= a\vartheta(\mathbf{X}+\bm{\Delta}\mathbf{X}^{(1)})
+b\vartheta(\mathbf{X}+\bm{\Delta}\mathbf{X}^{(2)})
\end{equation*}
and
\begin{equation*}
\lim\limits_{\|\bm{\Delta}\mathbf{X}\|\rightarrow0}\dfrac{\varsigma(\mathbf{X}+\bm{\Delta}\mathbf{X})}
{\|\bm{\Delta}\mathbf{X}\|_F}=0.
\end{equation*}
Then we define
\begin{equation*}
\nabla^2 h(\mathbf{X})\bm{\Delta}\mathbf{X}=\vartheta(\mathbf{X}+\bm{\Delta}\mathbf{X}).
\end{equation*}
If $\langle\nabla^2 h(\mathbf{X})\bm{\Delta}\mathbf{X},\bm{\Delta}\mathbf{X}\rangle>0$ for any $\bm{\Delta}\mathbf{X}\in\mathbb{Q}^{n_1\times n_2}$ and $\bm{\Delta}\mathbf{X}\neq\mathbf{0}$, then we call $\nabla^2 h$ is positive definite at $\mathbf{X}$.

\begin{prop}[\cite{qi2020quaternion}]\label{prop2.2}
Suppose that $h:\mathbb{Q}^{n_1\times n_2} \rightarrow\mathbb{R}$ is twice continuously differentiable. If $
\nabla h(\mathbf{X}^\sharp)=0$ and $\langle\nabla^2 h(\mathbf{X}^\sharp)\bm{\Delta}\mathbf{X},\bm{\Delta}\mathbf{X}\rangle>0$ for any $\bm{\Delta}\mathbf{X}\in\mathbb{Q}^{n_1\times n_2}$ and $\bm{\Delta}\mathbf{X}\neq\mathbf{0}$, then $\mathbf{X}^\sharp$ is a  minimizer value point of $h$.

\end{prop}


\section{The NRQMC model and algorithm}
\label{sec3}
In this section, we first clarify our motivation and propose a nonconvex robust quaternion matrix completion (NRQMC) model.  Then we establish an algorithm for solving the model and give the convergence analysis.

\subsection{The motivation}
\begin{itemize}
\item  \textbf{The low rank part:}
\end{itemize}
It is known that quaternion rank is the $L_0$-norm of singular value vector, while the  QNN is the $L_1$-norm of singular value vector.
We see from the quaternion singular value thresholding (QSVT) \cite[Theorem 2]{chen2020lowrank} that QSVT shrinks each singular value equally, i.e., each singular value is subtracted by the same threshold. This means that   QNN over-penalizes the large singular value.  Since QNN is a convex surrogate for quaternion matrix rank \cite{jia2019robust}, it is difficult to overcome the deficiency of QNN in quaternion rank approximation under the convex setting. Inspired by the recent extensive research on nonconvex techniques for robust  real tensor recovery, we perform the nonconvex MCP function to the singular value of quaternion matrix.
The MCP function is initially used for sparse estimation and variable
selections \cite{fan2001variable,zhang2010nearly}.  Next, Qiu et al. \cite{qiu2021nonlocal} adopted these two functions to construct  a tighter tensor rank approximation under transformed tensor
SVD framework for the robust tensor recovery problem.

The MCP function $\Phi_{c,\eta}$ is nonconvex and satisfies the following properties:

(1)  $\Phi_{c,\eta}$  is a mapping from $\mathbb{R}$ to $\mathbb{R}$ and $\Phi_{c,\eta}(0)=0$.

(2) $\Phi_{c,\eta}$  is proper, lower semicontinuous  and symmetric with respect to $y$-axis.

(3) $\Phi_{c,\eta}$  is concave and monotonically nondecreasing on $[0,+\infty)$.

Considering  MCP function $\Phi_{c,\eta}$ is symmetric with respect to $y$-axis, we now give an expression  of this function on $[0,+\infty)$
\begin{equation}\label{3.1}
\Phi_{c,\eta}(x)=c\int_0^x\max\{1-t/(c\eta),0\}dt=
\left\{
\begin{array}{ll}
	cx-\dfrac{1}{2\eta}x^2,&\text{if}~0\leq x\leq c\eta\\
	\dfrac{c^2\eta}{2},&\text{if}~x>c\eta,
\end{array}
\right.
\end{equation}
where $c,\eta>0$. The parameters $c$ and $\eta$ control   the
steepness of the quadratic function and the level of concavity. Besides, $c$ and $\eta$ influence the domain of $\Phi_{c,\eta}(x)$ to be a quadratic function or a constant function. When $c$ is fixed, the peak value of $\Phi_{c,\eta}(x)$ and the domain of the quadratic function are larger if $\eta$ is larger. This phenomenon also occurs
when $\eta$ is fixed and $c$ varies.

For a quaternion matrix $\mathbf{X}\in\mathbb{Q}^{n_1\times n_2}$, we define a nonconvex surrogate for the quaternion matrix rank by using the MCP function as follows
\begin{equation}\label{3.2}
\|\mathbf{X}\|_{\text{MCP}}=\sum\limits_{i=1}^{n_{(2)}}\Phi_{c,\eta}(\sigma_i(\mathbf{X})),~
c,\eta>0,
\end{equation}
where $\sigma_i(\mathbf{X})$ is the $i$-th singular value of $\mathbf{X}$ with $\sigma_1(\mathbf{X})\geq\cdots\geq\sigma_{n_{(2)}}(\mathbf{X})$.

\begin{figure}[!ht]
\begin{center}
	\subfloat[ ]{
		\includegraphics[width=6.0cm]{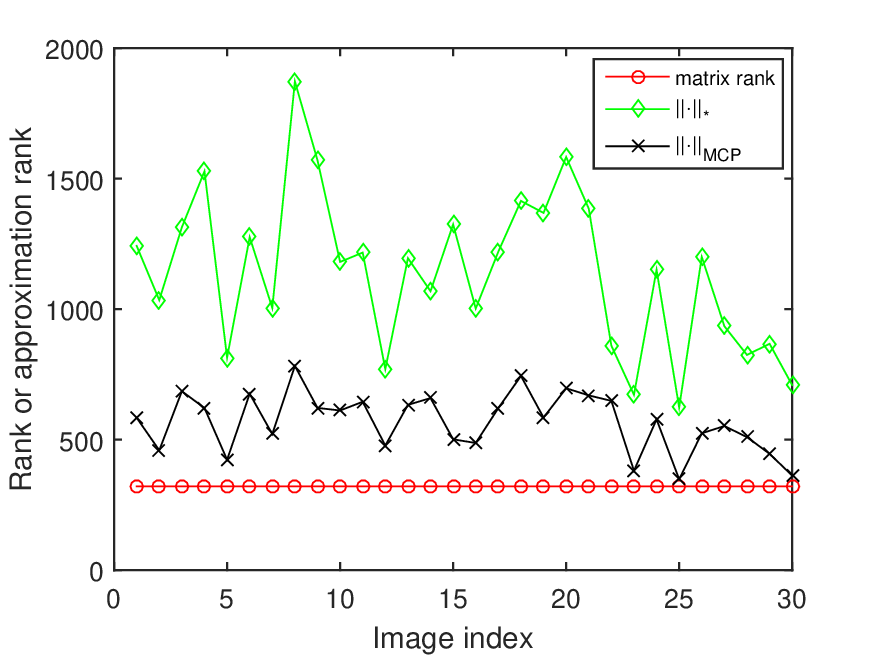}
	}\hspace{5mm}
	\subfloat[ ]{
		\includegraphics[width=6.0cm]{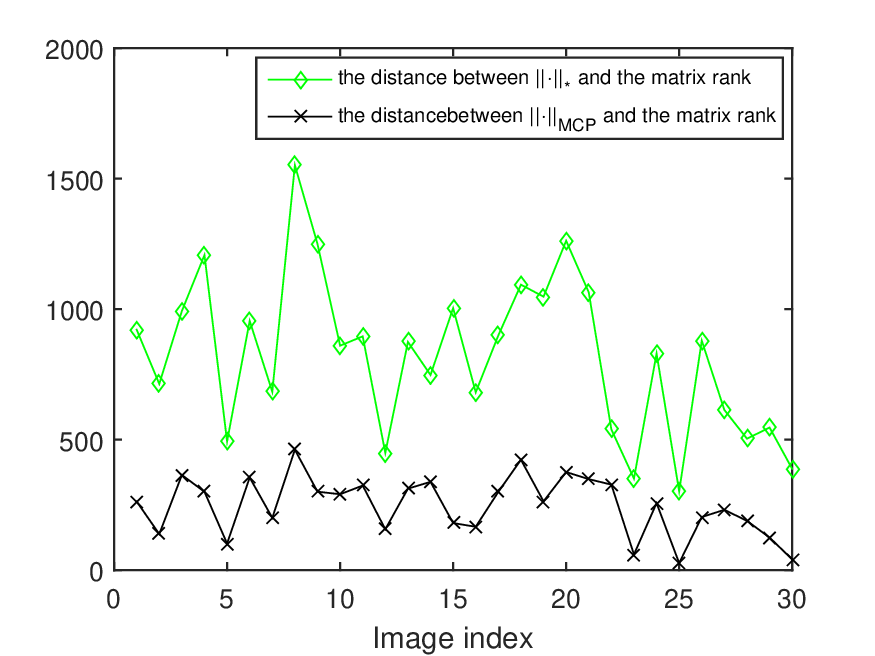}
	}
\end{center}
\caption{\label{convex_qnn_rank} Comparison of the quaternion matrix rank, the QNN and the MCP approximation rank function given by (\ref{3.2}).}
\end{figure}

Now we illustrate the reason for using the MCP function as the nonconvex quaternion matrix rank approximation. We randomly choose thirty color images from Berkeley Segmentation Dataset (BSD) \cite{martin2001database}, each with a size of $321\times 481$ or $481\times 321$. In Figure \ref{convex_qnn_rank} (a), for each image, we show the comparison of the quaternion matrix rank, the QNN  and the MCP approximation rank ($c=2,\eta=1.5$) given by (\ref{3.2}). We also show the distance between $\|\cdot\|_*$, $\|\cdot\|_{\text{MCP}}$ and the quaternion matrix rank for each image in Figure \ref{convex_qnn_rank} (b). We see from Figure \ref{convex_qnn_rank} that the result obtained by (\ref{3.2}) gives a tighter approximation to the quaternion matrix rank than the QNN for each image. The comparison results imply that our proposed nonconvex surrogate  (\ref{3.2}) appears to be better approximation for the quaternion matrix rank. The following numerical experiment results in Section \ref{sec5} show that the efficiency of MCP approximation rank function.

\begin{itemize}
\item  \textbf{The sparse part:}
\end{itemize}
As a convex relaxation of the quaternion $L_0$-norm, the quaternion $L_1$-norm has been used for sparsity in robust QMC model \cite{jia2019robust}. However, the obtained solution by the quaternion $L_1$-norm minimization may be
suboptimal to the original quaternion $L_0$-norm minimization since the quaternion $L_1$-norm is a coarse approximation
of the quaternion $L_0$-norm. It is shown that the use of nonconvex sparsity formulations can improve the realism of models and enhance their performance in
signal and image processing \cite{chartrand2007exact,xu2010regularization} and tensor robust PCA \cite{zhao2020nonconvex}. This inspires us to explore the nonconvex sparsity formulations in the robust QMC.

Considering the sparsity enhancement by using the quaternion $L_1$-norm minimization in
the robust QMC \cite{jia2019robust}, we employ the quaternion $L_p$-norm as the nonconvex relaxation of the quaternion $L_0$-norm.  For a quaternion matrix $\mathbf{X}\in\mathbb{Q}^{n_1\times n_2}$, the quaternion $L_p$-norm for $\mathbf{X}$ is specified by
\begin{equation}\label{3.3}
\|\mathbf{X}\|_p=\left(\sum\limits_{i=1}^{n_1}\sum\limits_{j=1}^{n_2}
|\mathbf{x}_{ij}|^p\right)^{\frac{1}{p}},~0<p<1.
\end{equation}
It is obvious to see that $\lim\limits_{p\rightarrow 1}\|\mathbf{X}\|_p=\|\mathbf{X}\|_1$.

\subsection{The proposed model and algorithm}
In order to overcome the disadvantages of QNN and quaternion $L_1$-norm, we  use two nonconvex replacements given by (\ref{3.2}) and (\ref{3.3}) to construct the following nonconvex robust QMC  (NRQMC) model
\begin{equation}\label{3.4}
\min\limits_{\mathbf{L},\mathbf{S}\in\mathbb{Q}^{n_1\times n_2}}\|\mathbf{L}\|_{\text{MCP}}+
\lambda \|\mathbf{S}\|_p^p,~~\text{s.t.}~~
\mathcal{P}_\Omega(\mathbf{L}+\mathbf{S})=\mathcal{P}_\Omega(\mathbf{X}).
\end{equation}
For convenience, we denote the model (\ref{3.4}) by ``MCP-L$_p$".
The model (\ref{3.4}) reduces to the following NRQPCA  problem when $\Omega$ is the entire set of indices
\begin{equation}\label{3.5}
\min\limits_{\mathbf{L},\mathbf{S}\in\mathbb{Q}^{n_1\times n_2}}\|\mathbf{L}\|_{\text{MCP}}+
\lambda \|\mathbf{S}\|_p^p,~~\text{s.t.}~~
\mathbf{L}+\mathbf{S}=\mathbf{X},
\end{equation}
and it  is equivalent to the following QMC problem when there is no corruption, i.e., $\mathbf{S}=\mathbf{0}$,
\begin{equation}\label{3.6}
\min\limits_{\mathbf{L}\in\mathbb{Q}^{n_1\times n_2}}\|\mathbf{L}\|_{\text{MCP}},~~\text{s.t.}~~
\mathcal{P}_\Omega(\mathbf{L})=\mathcal{P}_\Omega(\mathbf{X}).
\end{equation}

In the following, we develop the ADMM method to the nonconvex problem (\ref{3.4}).  ADMM decomposes a large global problem into a series of smaller subproblems, and coordinates the solutions of subproblems to compute the globally optimal solution. For more details, one can refer to \cite{candes2011robust}. It is noted that the quaternion $L_p$-norm in (\ref{3.4})  forces any entry of an optimal solution $\mathbf{S}$ in the
unobserved set $\Omega^\perp$ to be zero. Without loss of generality, we can assume that the unobserved data may be appropriate values such that $\mathcal{P}_{\Omega^\perp}(\mathbf{L}+\mathbf{S})=\mathcal{P}_{\Omega^\perp}(\mathbf{X})$. Then, the linear projection
operator constraint in (\ref{3.4}) is simply replaced by an equation $\mathbf{X}=\mathbf{L}+\mathbf{S}$. Thus, the problem (\ref{3.4}) can
be reformulated as
\begin{equation}\label{3.7}
\min\limits_{\mathbf{L},\mathbf{S}\in\mathbb{Q}^{n_1\times n_2}}\|\mathbf{L}\|_{\text{MCP}}+
\lambda \|\mathcal{P}_\Omega(\mathbf{S})\|_p^p,~~\text{s.t.}~~
\mathbf{L}+\mathbf{S}=\mathbf{X}.
\end{equation}
The augmented Lagrange
function of (\ref{3.7}) is given by
\begin{small}
\begin{eqnarray}
\!\!\!\!\hspace{-2mm}\mathscr{L}_\mu(\mathbf{L},\mathbf{S},\mathbf{M})
\!\!\!\!&=&\!\!\!\!
\|\mathbf{L}\|_{\text{MCP}}+
\lambda \|\mathcal{P}_\Omega(\mathbf{S})\|_p^p\!+\! \dfrac{\mu}{2}\Big\|\mathbf{L}\!+\!\mathbf{S}\!-\!\mathbf{X}\!+\!\frac{\mathbf{M}}{\mu}\Big\|_F^2
\!	-\!\dfrac{1}{2\mu}\|\mathbf{M}\|_F^2,\label{3.8}
\end{eqnarray}
\end{small}
where $\mathbf{M}$ is a Lagrange multiplier and  $\mu$ is the penalty parameter for linear constraints to
be satisfied.
Under the ADMM framework, $\mathbf{L}$, $\mathbf{S}$, and $\mathbf{M}$ can be alternately updated as follows:
\begin{footnotesize}
\begin{equation}\label{3.9}
\left\{
\begin{array}{l}
	\text{Step 1:}~~\mathscr{P}_\Omega(\mathbf{S}_{k+1})=\arg\underset{\mathbf{S}}{\min}\lambda \|\mathcal{P}_\Omega(\mathbf{S})\|_p^p+ \dfrac{\mu_k}{2}\|\mathcal{P}_\Omega(\mathbf{L}_{k}+\mathbf{S}-\mathbf{X}
	+\frac{\mathbf{M}_k}{\mu_k})\|_F^2\medskip\\
	\text{Step 2:}~~\mathscr{P}_{\Omega^\perp}(\mathbf{S}_{k+1}) =\arg\underset{\mathbf{S}}{\min}\dfrac{\mu_k}{2}\|\mathcal{P}_{\Omega^\perp}(\mathbf{L}_{k}+\mathbf{S}-\mathbf{X}
	+\frac{\mathbf{M}_k}{\mu_k})\|_F^2\medskip\\
	\text{Step 3:}~~\mathbf{L}_{k+1}=\arg\underset{\mathbf{L}}{\min}\|\mathbf{L}\|_{\text{MCP}}+
	\dfrac{\mu_k}{2}\|\mathbf{L}+\mathbf{S}_{k+1}-\mathbf{X}+\frac{\mathbf{M}_k}{\mu_k}\|_F^2\medskip\\
	\text{Step 4:}~~\mathbf{M}_{k+1}=\mathbf{M}_k+\mu_k(\mathbf{L}_{k+1}+\mathbf{S}_{k+1}-\mathbf{X}).
\end{array}
\right.
\end{equation}
\end{footnotesize}

\medskip
{\bf $\mathbf{S}$-subproblem:}  Before giving the solution of Step 1, we first extend the $L_p$-norm minimization problem from the real number field to the quaternion skew-field.

For a given quaternion $\mathbf{y}=y_0+y_1\mathbf{i}+y_2\mathbf{j}+y_3\mathbf{k}\in\mathbb{Q}$, let
\begin{equation}\label{3.14}
g(\mathbf{x})=\dfrac{1}{2}|\mathbf{x}-\mathbf{y}|^2+\nu|\mathbf{x}|^p,
\end{equation}
where $\nu>0$, $0<p<1$,  and $\mathbf{x}=x_0+x_1\mathbf{i}+x_2\mathbf{j}+x_3\mathbf{k}\in\mathbb{Q}$.

The optimal solution of the minimization problem {\rm(\ref{3.14})} is given by
the following quaternion generalized soft-thresholding {\rm(QGST)} operator.
\begin{thm}\label{thm3.2}
For the given $\nu> 0$ and  $p$ $(0<p<1)$, an optimal solution of the minimization problem {\rm(\ref{3.14})} is given by
the quaternion generalized soft-thresholding {\rm(QGST)} operator, which is defined by
\begin{equation}\label{3.22}
	\mathrm{QGST}(\mathbf{y},\nu,p)=\left\{
	\begin{array}{ll}
		0,&\text{if}~|\mathbf{y}|\leq\tau_p^{\mathrm{QGST}}(\nu),\\
		\mathrm{signQ}(\mathbf{y})T_p^{\mathrm{QGST}}(|\mathbf{y}|,\nu),
		&\text{if}~|\mathbf{y}|>\tau_p^{\mathrm{QGST}}(\nu),\\
	\end{array}
	\right.
\end{equation}
where $\tau_p^{\mathrm{QGST}}(\nu)=(2\nu (1-p))^{\frac{1}{2-p}}+\nu p(2\nu (1-p))^{\frac{p-1}{2-p}}$ is a threshold value, and $T_p^{\mathrm{QGST}}(|\mathbf{y}|,\nu)$ can be obtained by solving the following equation
\begin{equation}\label{3.23}
	T_p^{\mathrm{QGST}}(|\mathbf{y}|,\nu)-|\mathbf{y}|+\nu p(T_p^{\mathrm{QGST}}(|\mathbf{y}|,\nu))^{p-1}=0.
\end{equation}
\end{thm}

To solve (\ref{3.14}), we propose an iterative algorithm, which is summarized as in Algorithm \ref{3.1}.

\begin{algorithm}[!ht]
\caption{Quaternion Generalized Soft-Thresholding (QGST)}\label{alg3.2}
{\bf Input:} $\mathbf{y}$, $\lambda$, $p$, $J=2$ or $3$.
\begin{algorithmic}[1]
	\STATE $\tau_p^{\mathrm{QGST}}(\lambda)=(2\lambda (1-p))^{\frac{1}{2-p}}+\lambda p(2\lambda (1-p))^{\frac{p-1}{2-p}}$;
	\STATE \textbf{if}~$|\mathbf{y}|<\tau_p^{\mathrm{QGST}}(\lambda)$
	\STATE \quad $\mathrm{QGST}(\mathbf{y},\lambda,p)=0$;
	\STATE \textbf{else}
	\STATE \quad $k=0$, $t^{(k)}=|\mathbf{y}|$
	\STATE \quad \textbf{for}~$k=0,1,\cdots, J$~\textbf{do}
	\STATE \quad\quad  $t^{(k+1)}=|\mathbf{y}|-\lambda p(t^{(k)})^{p-1}$;
	\STATE \quad\quad $k=k+1$;
	\STATE \quad \textbf{end~for}
	\STATE \quad $\mathrm{QGST}(\mathbf{y},\lambda,p)=\mathrm{signQ}(\mathbf{y})t^{(k)}$;
	\STATE \textbf{end~if}
\end{algorithmic}
{\bf Output:} $\mathrm{QGST}(\mathbf{y},\lambda,p)$.
\end{algorithm}

\begin{rem}\label{rem3.1}
When $p=1$, {\rm QGST} will converge after one iteration. Since
\begin{equation*}
	\lim\limits_{p\rightarrow 1}\tau_p^{\mathrm{QGST}}(\nu)= \lim\limits_{p\rightarrow 1}
	\nu p(2\nu (1-p))^{\frac{p-1}{2-p}}=\nu\lim\limits_{p\rightarrow 1}(2\nu )^{p-1}
	\lim\limits_{p\rightarrow 1}(1-p )^{p-1}=\nu.
\end{equation*}
the thresholding value of {\rm QGST} will become $\nu$, and the {\rm QGST}
function becomes
\begin{equation}\label{3.24}
	\mathrm{QGST}(\mathbf{y},\nu,1)=\left\{
	\begin{array}{ll}
		0,&\text{if}~|\mathbf{y}|\leq\nu,\\
		\mathrm{signQ}(\mathbf{y})(|\mathbf{y}|-\nu),
		&\text{if}~|\mathbf{y}|>\nu.
	\end{array}
	\right.
\end{equation}
One can see that {\rm(\ref{3.24})} is just the soft-thresholding function proposed by Jia et al. {\rm\cite{jia2022nonlocal}}.
\end{rem}

According to Step 1, it follows that
\begin{equation}\label{3.25}
(\mathbf{S}_{k+1})_{ij}=\arg\min \frac{1}{2}|\mathbf{S}_{ij}-(\mathbf{Z}_k)_{ij}|^2+\frac{\lambda}{\mu_k} |(\mathbf{S})_{ij}|^p,~(i,j)\in\Omega,
\end{equation}
where $\mathbf{Z}_k=\mathbf{X}-\mathbf{L}_{k}
-\frac{\mathbf{M}_k}{\mu_k}$.
According to Theorem \ref{thm3.2} and Step 2, the solution with respect to $\mathbf{S}_{k+1}$ is given by
\begin{equation}\label{3.26}
\mathbf{S}_{k+1} = \mathcal{P}_\Omega ( \mathrm{QGST}( (\mathbf{Z}_k),\dfrac{\lambda}{\mu_k},p) ) + \mathcal{P}_{\Omega^\perp} (\mathbf{Z}_k).
\end{equation}

\medskip
{\bf $\mathbf{L}$-subproblem:}  According to the relation (\ref{3.9}), we need to solve the following subproblem:
\begin{equation}\label{3.10}
\min\dfrac{1}{\mu_k}\|\mathbf{L}\|_{\text{MCP}}+
\dfrac{1}{2}\|\mathbf{L}-\mathbf{Y}_k\|_F^2,
\end{equation}
where $\mathbf{Y}_k=\mathbf{X}-\mathbf{S}_{k+1}-\frac{\mathbf{M}_k}{\mu_k}$.

The solution of the subproblem (\ref{3.10}) is given by the following theorem.


\begin{thm}[Quaternion thresholding operator
for    MCP penalty]\label{thm3.1}
Let $\mathbf{Y}_k=\mathbf{U}_k\mathbf{\Sigma}_k\mathbf{V}_k^*$ be the {\rm QSVD} of $\mathbf{Y}_k\in\mathbb{Q}^{n_1\times n_2}$. Then an optimal solution of the  minimization problem {\rm(\ref{3.10})} is given by
$
	\mathbf{L}_{\clubsuit}=\mathbf{U}_k \mathbf{\Sigma}_{\frac{1}{\mu_k}\Phi_{c,\eta}}\mathbf{V}_k^*$,
where
\begin{equation}\label{new1}
	\mathbf{\Sigma}_{\frac{1}{\mu_k}\Phi_{c,\eta}}=\mathrm{Diag} \Big(
	\mathrm{Prox}_{\frac{1}{\mu_k}\Phi_{c,\eta}}\big((\mathbf{\Sigma}_k)_{1,1}\big),\cdots,
	\mathrm{Prox}_{\frac{1}{\mu_k}\Phi_{c,\eta}}\big((\mathbf{\Sigma}_k)_{n_{(2)},n_{(2)}}\big)\Big)
\end{equation}
with
$(\mathbf{\Sigma}_k)_{i,i}$ being the $(i,i)$th entry of $\mathbf{\Sigma}_k$, and $\mathrm{Prox}_{\frac{1}{\mu_k}\Phi_{c,\eta}}\big((\mathbf{\Sigma}_k)_{i,i}\big)$ being the proximal mapping associated with $\frac{1}{\mu_k}\Phi_{c,\eta}$ at $(\mathbf{\Sigma}_k)_{i,i}$.
\end{thm}

Next we give the proximal mapping of $\mu \Phi_{c,\eta}$ in (\ref{3.1}). For $0<\mu<\eta$,  $\mathrm{Prox}_{\mu\Phi_{c,\eta}}(y)$ is given by
\begin{equation}\label{3.13}
\mathrm{Prox}_{\mu\Phi_{c,\eta}}(y)
=\left\{
\begin{array}{ll}
	0,&\text{if}~|y|\leq c\mu,\\
	\dfrac{\mathrm{sign}(y)(|y|-c\mu)}{1-\mu/\eta},&\text{if}~c\mu<|y|\leq c\eta,\\
	y,&\text{if}~|y|>c\eta,
\end{array}
\right.
\end{equation}
where $\mathrm{sign}(y)$ equals to $1$, $0$, and $-1$ if $y>0$, $y=0$, and $y<0$, respectively.
It is noted that $\eta>\mu$ and $c > 0$, where  $\eta>\mu$ is to guarantee the meaningful solution in the proximal mapping about
the MCP function and $c>0$ is to guarantee the concavity of the quadratic function
on $[0,c\eta]$.

For the sake of clarity, we summarize the update of $\mathbf{L}$ in Algorithm \ref{alg3.1}. 

\begin{algorithm}[!ht]
\caption{Update $\mathbf{L}$}\label{alg3.1}
{\bf Input:} quaternion matrices $\mathbf{S}_k$ and $\mathbf{M}_k\in\mathbb{Q}^{n_1\times n_2}$  and scalars $\mu_k$, $c$, $\eta$.

\begin{algorithmic}[1]
	\STATE $\mathbf{Y}_k=\mathbf{X}-\mathbf{S}_{k+1}-\frac{\mathbf{M}_k}{\mu_k}$;
	\STATE $[\mathbf{U}_k~\mathbf{\Sigma}_k~\mathbf{V}_k]=\mathrm{QSVD}(\mathbf{Y}_k)$;
	\STATE Compute $\mathbf{\Sigma}_{\frac{1}{\mu_k}\Phi_{c,\eta}}$ via  (\ref{new1});
	\STATE $\mathbf{L}_{k+1}=\mathbf{U}_k \mathbf{\Sigma}_{\frac{1}{\mu_k}\Phi_{c,\eta}}\mathbf{V}_k^*$;
\end{algorithmic}
{\bf Output:} $\mathbf{L}_{k+1}$.
\end{algorithm}

The whole ADMM algorithm for solving model (\ref{3.4}) is summarized in Algorithm \ref{alg3.3}.

\begin{algorithm}[!ht]
\caption{ADMM for Solving NRQMC Model}\label{alg3.3}
{\bf Input:} The observed quaternion matrix $\mathbf{X}\in\mathbb{Q}^{n_1\times n_2}$ with $\Omega$ (the index of observed entries); balanced parameter $\lambda$; parameters $p$, $c$ and $\eta$.\\
{\bf Initialize:} $\mathbf{L}_0=\mathbf{S}_0=\mathbf{M}_0=0$, $\mu_0=10^{-4}$,
\begin{algorithmic}[1]
	\STATE \textbf{while} not converged and $k<500$ \textbf{do}
	\STATE \quad Update $\mathbf{S}_{k+1}$ via (\ref{3.26});
	\STATE \quad Update $\mathbf{L}_{k+1}$ according to Algorithm (\ref{alg3.1});
	\STATE \quad Update $\mathbf{M}_{k+1}$ via (\ref{3.9});
	\STATE \quad Update $\mu_{k+1}$ by $\mu_{k+1}=\min\{1.2*\mu_k,10^8\}$;
	\STATE \quad Check the convergence
	\begin{equation}\label{3.28}
		\max\{\|\mathbf{L}_{k+1}-\mathbf{L}_k\|_F,\|\mathbf{S}_{k+1}-\mathbf{S}_k\|_F, \|\mathbf{L}_{k+1}+\mathbf{S}_{k+1}-\mathbf{X}\|_F\}\leq tol;
	\end{equation}
	\STATE \quad Set $k=k+1$;
	\STATE \textbf{end~while}
\end{algorithmic}
\textbf{Output:} $\mathbf{L}$ and $\mathbf{S}$.
\end{algorithm}
%
The computation cost of Algorithm \ref{alg3.3} mainly lies in updating $\mathbf{L}_{k+1}$, in which the QSVD of quaternion matrix $\mathbf{Y}_k$ is needed. It is well known that  the computational complexity of the QSVD for an $n_1\times n_2$ quaternion matrix is $O(n_1n_2\min \{n_1, n_2\})$.  Therefore, at each iteration, the total computational complexity of the proposed algorithm is $O(n_1n_2\min \{n_1, n_2\})$. For simplicity, if $n_1=n_2=n$, the cost at each iteration is $O(n^3)$.


\subsection{Convergence analysis}
Now we give the convergence analysis of Algorithm \ref{alg3.3}. 

\begin{thm}\label{thm4.1}
Let the sequence $\mathbf{W}_k = \{\mathbf{L}_k, \mathbf{S}_k,\mathbf{M}_k\}$ be generated by {\rm Algorithm \ref{alg3.3}}. Then the
accumulation point $\mathbf{W}_\clubsuit = \{\mathbf{L}_\clubsuit, \mathbf{S}_\clubsuit,\mathbf{M}_\clubsuit\}$ of $\mathbf{W}_k$ is a {\rm KKT} stationary point, i.e., $\mathbf{W}_\clubsuit$ satisfies
the following {\rm KKT} conditions:
\begin{eqnarray*}
	&& \nabla \|\mathbf{L}\|_{\mathrm{MCP}}|_{\mathbf{L}=\mathbf{L}_\clubsuit}+\mathbf{M}_\clubsuit=0,
	\langle\mathcal{P}_{\Omega}(\mathbf{S}_\clubsuit),\mathcal{P}_{\Omega}(\mathbf{M}_\clubsuit)\rangle
	+\lambda p\|\mathcal{P}_{\Omega}(\mathbf{S}_\clubsuit)\|_p^p=0,\\
	&&\mathcal{P}_{\Omega^\perp}(\mathbf{M}_\clubsuit)=0,
	\mathrm{L}_\clubsuit+\mathbf{S}_\clubsuit-
	\mathbf{X}=0.
\end{eqnarray*}
\end{thm}

The proof of  Theorem \ref{thm4.1}   can be found in the Appendix.

\section{NSS-based method for large-scale color image and video inpainting}
\label{sec4}
In this section, we develop a robust quaternion completion method that jointly exploits low-rankness and nonlocal self-similarity (NSS).
Given an observed quaternion tensor $\mathcal{X}\in\mathbb{Q}^{n_1 \times n_2 \times n_3}$ and the set of the indices of known elements $\Omega$, we obtain an estimated image $\widetilde{\mathcal{L}}$ by a certain robust quaternion completion method.
For each slice $\widetilde{\mathcal{L}}(:,:,j) \in \mathbb{Q}^{n_1 \times n_2}$, we divide it into $M$ subblocks with size $s \times s$ and overlap $l$, where $M=\left\lceil{(n_1-l)} / {(s-l)}\right\rceil \times \left\lceil {(n_2-l)} / {(s-l)}\right\rceil$. We therefore obtain $Mn_3$ subblocks and then apply the K-means++ algorithm to cluster these subblocks into $N$ groups $\{ \widetilde{\mathbf{L}}^i\}_{i=1}^N$, where $\widetilde{\mathbf{L}}^i \in \mathbb{Q}^{s^2 \times m_i}$, each column of $\widetilde{\mathbf{L}}^i$ is the vectorization of the corresponding subblock, $m_i$ is the number of subblocks in the $i$th group, and $\sum_{i=1}^{N} m_i =M$. According to the coordinates of subblocks $\widetilde{\mathbf{L}}^i$, the sub-data $\mathbf{X}^i$ and the index set $\Omega ^i$ are generated from $\mathcal{X}$ and $\Omega$, respectively. For each $i=1,\cdots,N$, we solve the following subproblem
\begin{equation}\label{sub-pro}
\min\limits_{\mathbf{L}^i,\mathbf{S}^i}\|\mathbf{L}^i\|_{\text{MCP}}+
\lambda \|\mathcal{P}_{\Omega^i}(\mathbf{S}^i)\|_p^p,~~\text{s.t.}~~
\mathbf{L}^i+\mathbf{S}^i=\mathbf{X}^i.
\end{equation}
Finally, we reshape each column of each $\mathbf{L}^i$ as a quaternion matrix of size $s \times s$ and then aggregate these subblocks to reconstruct the final data. In addition, the pixels on the positions where the subblocks overlap are set to be the average of the corresponding values. 

Now we propose the method combining NRQMC and NSS for robust quaternion completion as stated in Algorithm \ref{alg4.1}. Note that the third dimension of a quaternion matrix data is 1 and the number of subblocks is $M$. We use NRQMC-NSS2D and NRQMC-NSS3D to represent Algorithm \ref{alg4.1} for processing 2D and 3D quaternion data, respectively.
\begin{algorithm}[H]
\caption{NRQMC-NSS Algorithm for Robust Quaternion Completion}\label{alg4.1}
{\bf Step 1.} Input the observed data $\mathcal{X}$ and the index of observed entries $\Omega$, the initial estimatior $\widetilde{\mathcal{L}}$, the patch size $s$, the overlap size $l$.

{\bf Step 2.} Divide $\widetilde{\mathcal{L}}$ into $Mn_3$ subblocks and obtain $N$ clusters $\{ \widetilde{\mathbf{L}}^i \}_{i=1}^{N}$ using K-means++ algorithm. Obtain $\mathbf{X}^i$ and $\Omega ^i$ according to the same coordinates of subblocks in $\widetilde{\mathbf{L}}^i$.


{\bf Step 3.} Apply Algorithm \ref{alg3.3} to obtain $\mathbf{L}^i$ according to $\mathbf{X}^i$ and $\Omega^i$.

{\bf Step 4.} Reconstruct the final data $\mathcal{L}$ by aggregating all subblocks in $\mathbf{L}^{i}$.

\end{algorithm}

\section{Numerical experiments}
\label{sec5}
In this section, numerical experiments are conducted on color images and videos to demonstrate the effectiveness of our proposed model and algorithm.

The sampling ratio $SR$ indicates the percentage of observed entries. It is defined as $SR = |\Omega|/(n_1n_2)$, where the size of the image is of $n_1 \times n_2$, $\Omega$ is the set of randomly generated observations and $|\Omega|$ denotes its cardinality.
The corrupted observation is generated as follows.   First, a standard uniform noise is independently and randomly added into $\ell$ pixel locations of red, green, and blue channels
of color images for producing the noisy data. The sparsity of noise components is denoted
by $\gamma= \ell/(n_1n_2)$. 
Second, a percentage $SR$ is chosen randomly from the noisy data as the observation.

The intensity range of the real-world data is scaled into $[0,1]$. All experiments are performed in the MATLAB 2024a environment and run on a desktop computer (Intel Core i7-11700, @ 2.5GHz, 32G RAM).

\subsection{Parameters Setting}
For the proposed models, six parameters are used to control the performance, i.e., $\lambda$, $\eta$, $c$, $p$, $s$, and $l$.
The regularization parameter $\lambda$ is used to balance the low-rank and sparse parts and is set to be $\lambda = 1/\sqrt{SR\cdot \max(n_1,n_2)}$. 

\begin{figure}[H]
\centering
\subfloat[\footnotesize Parameters $\eta$ and $c$]{\label{subfig:para_nc}
	\includegraphics[width=0.25\linewidth]{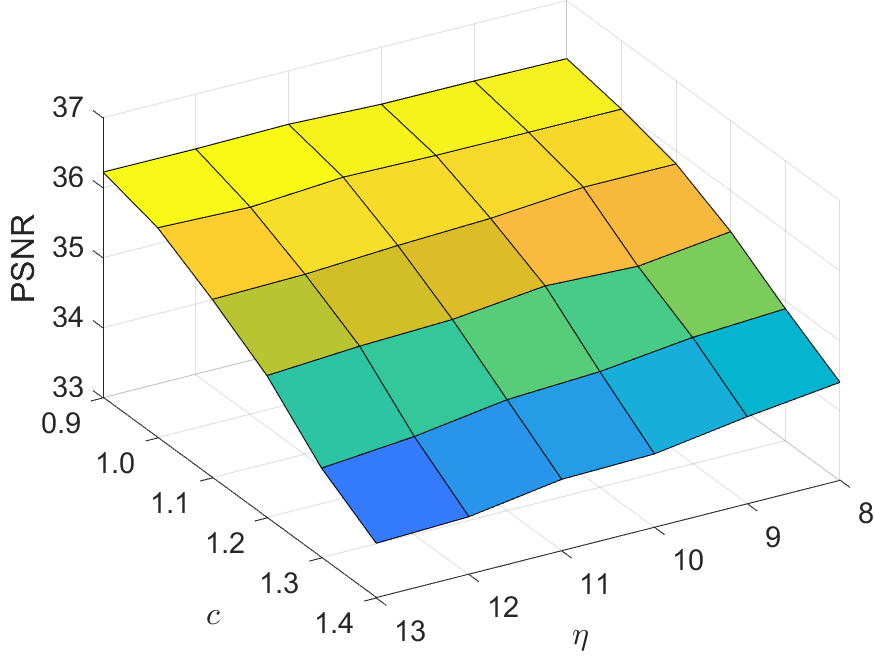}
}\hspace{10pt}
\subfloat[\footnotesize Parameter $p$]{\label{subfig:para_p}
	\includegraphics[width=0.25\linewidth]{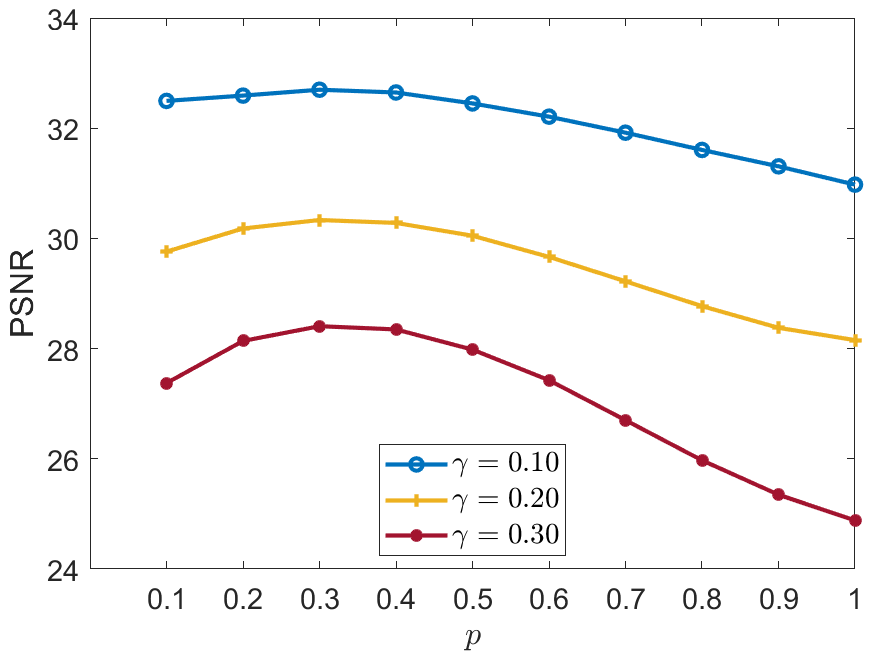}
}\hspace{10pt}
\subfloat[\footnotesize Size of subblocks $s$]{\label{subfig:para_w_sz}
	\includegraphics[width=0.25\linewidth]{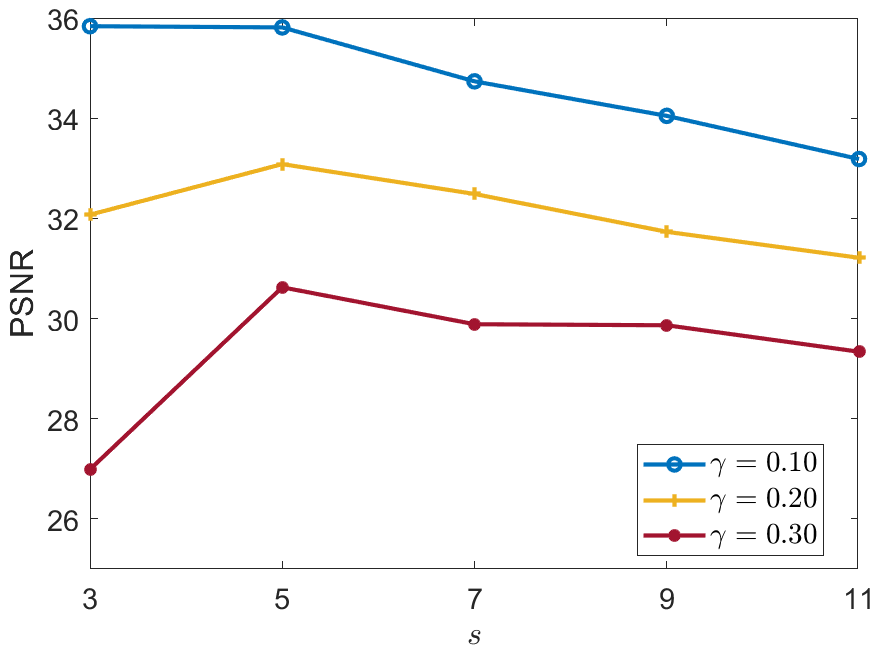}
}
\caption{Performance of the proposed methods with different parameters on \textit{Lena} under $SR=0.2$. }\label{fig:para_set}
\end{figure}

We take \textit{Lena} as an example to evaluate the influence of the parameters $p$, $\eta$, $c$, and the size of search window $s$, where $SR=0.2$.
Figure \ref{fig:para_set}\subref{subfig:para_nc} shows the PSNR values with respect to the parameters $\eta$ and $c$, where $\gamma$ is set to be 0.1. Parameter $\eta$ is chosen from $\left\lbrace8,9,10,11,12,13\right\rbrace$ and $c$ is chosen from $\left\lbrace0.9,1.0,1.1,1.2,1.3,1.4\right\rbrace$. We can see that the influence of $\eta$ is small and NRQMC performs well when $c$ is small. Therefore, we set $\eta=13$ and $c=0.9$ in all experiments.
Figure \ref{fig:para_set}\subref{subfig:para_p} shows the performance of NRQMC with different $p$ under $\gamma=0.1,0.2,0.3$. It can be seen that NRQMC performs well when $p$ falls within the range of $[0.2, 0.4]$, and obtains the best PSNR value when $p$ is set to be 0.3. Therefore, $p$ is set to be 0.3 in all experiments.
The overlap $l$ is always set to be 1 in our experiments.
Figure \ref{fig:para_set}\subref{subfig:para_w_sz} shows the performance of NRQMC-NSS with different $s$ under $\gamma=0.1,0.2,0.3$. We can see that NRQMC-NSS performs best when $s=5$ in the case $\gamma=0.2,0.3$. Although the PSNR value of $s=3$ is higher than that of $s=5$, the improvement is very small. Therefore, $s$ is set to be 5 in all experiments.

\subsection{Color Images}
In this subsection, we select six color images to demonstrate the performance of the proposed methods, i.e., \textit{House} ($256 \times 256$), \textit{Peppers} ($512 \times 512$), \textit{Lena} ($512 \times 512$)\footnote{\href{http://www.eecs.qmul.ac.uk/~phao/IP/Images/}{http://www.eecs.qmul.ac.uk/~phao/IP/Images/}}, \textit{Kodak image2} ($512 \times 768$), \textit{Kodak image3} ($512 \times 768$), \textit{Kodak image12} ($512 \times 768$)\footnote{\href{https://r0k.us/graphics/kodak/}{https://r0k.us/graphics/kodak/}}. The $SR$ is set to be 0.5, 0.6, and 0.8 and the $\gamma$ is set to be 0.1 and 0.2.
We compare the proposed NRQMC and NRQMC-NSS2D with
practical quaternion matrix completion algorithm (PQMC) \cite{jia2022nonlocal},
tensor completion algorithm (TNN) \cite{song2020robust},
low-rank quaternion approximation with Geman function  (LRQA-G) \cite{chen2020lowrank},
and robust quaternion matrix completion with convex surrogates (QNN) \cite{jia2019robust}.
In addition, we incorporated PQMC-NSS2D to compare the capability of PQMC and NRQMC in solving subblock completion problems, which uses PQMC to solve the subproblems in Step 3 of Algorithm \ref{alg4.1}.
The recovered image by NRQMC is used as the initial guess of PQMC-NSS2D and NRQMC-NSS2D.
We comment here that the LRQA-G algorithm \cite{chen2020lowrank} separately designed the corresponding framework for image inpainting and image denoising. For fairness of comparison, we make some modifications by adding the quaternion $L_1$-norm characterization of the sparse constraint and using the ADMM framework as in Algorithm \ref{alg3.1}.
The parameters in the compared methods are manually adjusted in all experiments to obtain the best result according to the authors' suggestions.
To evaluate the performance of the proposed model and algorithm, except for visual quality, we adopt two quantitative assessment indices including the peak signal-to-noise ratio (PSNR), and the structural similarity index (SSIM).

\begin{figure}[H]
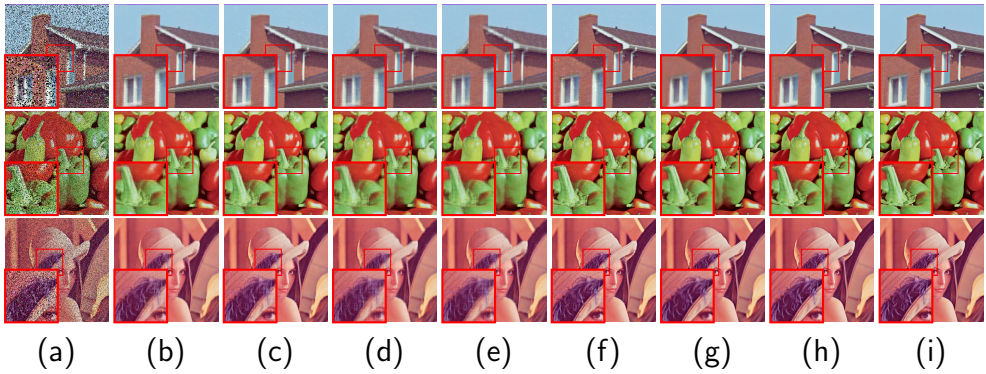

\centering
\subfloat{

}
\caption{Recovery results by different methods on \textit{House}, \textit{Peppers}, and \textit{Lena} (from top to bottle) under $SR=0.8$ and $\gamma = 0.2$. (a) Observed. (b) TNN. (c) LRQA-G. (d) QNN. (e) PQMC. (f) NRQMC. (g) PQMC-NSS2D. (h) NRQMC-NSS2D. (i) Original. }\label{fig:result_rgb1}%
\end{figure}

\begin{figure}[H]
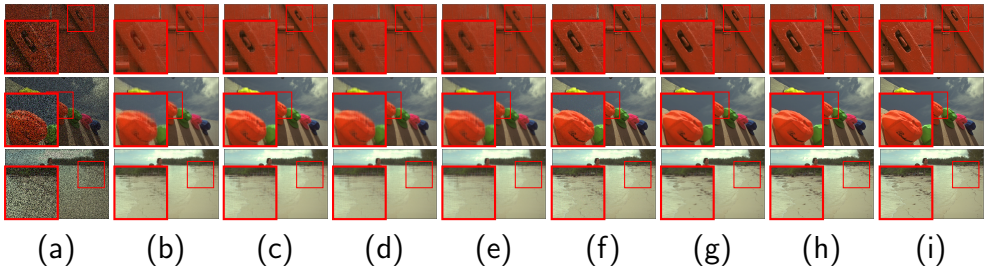

\centering
\subfloat{
	%
}
\caption{Recovery results by different methods on \textit{Kodak image12} under $SR=0.5$ and $\gamma = 0.1$. (a) Observed. (b) TNN. (c) LRQA-G. (d) QNN. (e) PQMC. (f) NRQMC. (g) PQMC-NSS2D. (h) NRQMC-NSS2D. (i) Original. }\label{fig:result_rgb2}%
\end{figure}

\begin{table}[H]
\centering
\caption{The PSNR and SSIM values of all methods for color images.}\label{tab:result_color_images}%
\scalebox{.52}{\setlength{\tabcolsep}{1.3mm}{\begin{tabular}{llllllllllllllllll}
			\hline
			&       &       & PSNR  &       &       &       &       &       &       &       & SSIM  &       &       &       &       &       &  \\
			\cline{4-10}\cline{12-18}          & SR    & $\gamma$ & TNN & LRQA & QNN   & PQMC  & NRQMC & PQMC & NRQMC &       & TNN & LRQA & QNN   & PQMC  & NRQMC & PQMC & NRQMC \\
			&       &       &       & -G     &      &       &       & -NSS2D   & -NSS2D   &       &       & -G     &    &       &       & -NSS2D   & -NSS2D \\
			\cline{1-10}\cline{12-18}    \textit{House} & 0.5   & 0.1   & 28.38  & 28.63  & 26.76  & 26.75  & 29.78  & 32.06  & \textbf{34.00} &       & 0.943  & 0.948  & 0.936  & 0.936  & 0.954  & 0.974  & \textbf{0.983} \\
			&       & 0.2   & 26.75  & 15.69  & 25.52  & 25.50  & 27.13  & 31.23  & \textbf{32.65} &       & 0.917  & 0.472  & 0.910  & 0.909  & 0.919  & 0.968  & \textbf{0.977} \\
			& 0.6   & 0.1   & 28.92  & 29.29  & 26.83  & 26.83  & 31.22  & 32.49  & \textbf{34.81} &       & 0.950  & 0.956  & 0.940  & 0.940  & 0.965  & 0.976  & \textbf{0.986} \\
			&       & 0.2   & 27.83  & 24.68  & 26.35  & 26.36  & 28.94  & 31.87  & \textbf{33.36} &       & 0.933  & 0.868  & 0.924  & 0.924  & 0.943  & 0.972  & \textbf{0.980} \\
			& 0.8   & 0.1   & 30.21  & 30.29  & 27.83  & 27.83  & 33.36  & 33.41  & \textbf{37.17} &       & 0.964  & 0.968  & 0.953  & 0.953  & 0.980  & 0.980  & \textbf{0.991} \\
			&       & 0.2   & 29.08  & 28.95  & 27.11  & 27.12  & 30.69  & 33.09  & \textbf{34.81} &       & 0.950  & 0.952  & 0.940  & 0.940  & 0.964  & 0.978  & \textbf{0.987} \\
			\cline{1-10}\cline{12-18}    \textit{Peppers} & 0.5   & 0.1   & 27.12  & 28.56  & 27.25  & 27.24  & 29.04  & 31.24  & \textbf{32.88} &       & 0.962  & 0.971  & 0.964  & 0.964  & 0.973  & 0.984  & \textbf{0.989} \\
			&       & 0.2   & 25.83  & 13.34  & 25.78  & 25.73  & 27.05  & 30.68  & \textbf{31.82} &       & 0.948  & 0.542  & 0.947  & 0.947  & 0.957  & 0.981  & \textbf{0.986} \\
			& 0.6   & 0.1   & 27.51  & 29.28  & 27.69  & 27.68  & 29.90  & 31.75  & \textbf{33.77} &       & 0.966  & 0.976  & 0.968  & 0.968  & 0.978  & 0.985  & \textbf{0.991} \\
			&       & 0.2   & 26.33  & 14.60  & 26.31  & 26.28  & 27.95  & 31.22  & \textbf{32.42} &       & 0.955  & 0.610  & 0.955  & 0.955  & 0.965  & 0.983  & \textbf{0.987} \\
			& 0.8   & 0.1   & 28.37  & 30.50  & 28.54  & 28.54  & 31.69  & 32.59  & \textbf{35.55} &       & 0.972  & 0.982  & 0.974  & 0.974  & 0.985  & 0.988  & \textbf{0.994} \\
			&       & 0.2   & 27.52  & 28.65  & 27.54  & 27.53  & 29.75  & 32.05  & \textbf{33.94} &       & 0.966  & 0.972  & 0.967  & 0.966  & 0.977  & 0.986  & \textbf{0.991} \\
			\cline{1-10}\cline{12-18}    \textit{Lena} & 0.5   & 0.1   & 28.22  & 28.75  & 27.32  & 27.30  & 29.30  & 31.91  & \textbf{33.56} &       & 0.969  & 0.971  & 0.966  & 0.966  & 0.972  & 0.985  & \textbf{0.990} \\
			&       & 0.2   & 26.70  & 13.51  & 26.01  & 25.98  & 27.36  & 31.29  & \textbf{32.36} &       & 0.957  & 0.518  & 0.953  & 0.953  & 0.958  & 0.982  & \textbf{0.986} \\
			& 0.6   & 0.1   & 28.62  & 29.33  & 27.51  & 27.51  & 30.28  & 32.23  & \textbf{34.54} &       & 0.972  & 0.975  & 0.968  & 0.968  & 0.978  & 0.986  & \textbf{0.991} \\
			&       & 0.2   & 27.45  & 15.04  & 26.57  & 26.55  & 28.39  & 31.95  & \textbf{33.30} &       & 0.963  & 0.604  & 0.960  & 0.960  & 0.966  & 0.985  & \textbf{0.989} \\
			& 0.8   & 0.1   & 29.88  & 30.62  & 28.46  & 28.45  & 32.61  & 33.54  & \textbf{36.88} &       & 0.979  & 0.982  & 0.974  & 0.974  & 0.987  & 0.989  & \textbf{0.995} \\
			&       & 0.2   & 28.74  & 28.90  & 27.48  & 27.47  & 30.38  & 32.99  & \textbf{34.97} &       & 0.973  & 0.972  & 0.968  & 0.968  & 0.978  & 0.988  & \textbf{0.992} \\
			\cline{1-10}\cline{12-18}    \textit{Kodak} & 0.5   & 0.1   & 29.39  & 30.30  & 28.93  & 28.92  & 30.60  & 31.59  & \textbf{32.92} &       & 0.981  & 0.985  & 0.980  & 0.980  & 0.987  & 0.989  & \textbf{0.992} \\
			\textit{image2} &       & 0.2   & 28.63  & 24.80  & 28.46  & 28.46  & 29.59  & 31.36  & \textbf{32.40} &       & 0.977  & 0.947  & 0.978  & 0.978  & 0.983  & 0.988  & \textbf{0.991} \\
			& 0.6   & 0.1   & 29.68  & 30.45  & 29.06  & 29.06  & 31.16  & 31.86  & \textbf{33.86} &       & 0.982  & 0.985  & 0.981  & 0.981  & 0.989  & 0.990  & \textbf{0.994} \\
			&       & 0.2   & 28.97  & 29.89  & 28.61  & 28.60  & 30.02  & 31.77  & \textbf{32.96} &       & 0.979  & 0.983  & 0.979  & 0.979  & 0.985  & 0.989  & \textbf{0.992} \\
			& 0.8   & 0.1   & 30.34  & 30.87  & 29.36  & 29.35  & 32.50  & 32.81  & \textbf{35.74} &       & 0.985  & 0.987  & 0.982  & 0.982  & 0.992  & 0.992  & \textbf{0.996} \\
			&       & 0.2   & 29.66  & 30.38  & 29.02  & 29.02  & 31.43  & 32.55  & \textbf{34.48} &       & 0.982  & 0.985  & 0.980  & 0.980  & 0.989  & 0.991  & \textbf{0.995} \\
			\cline{1-10}\cline{12-18}    \textit{Kodak} & 0.5   & 0.1   & 28.55  & 29.76  & 28.30  & 28.30  & 30.79  & 32.71  & \textbf{34.88} &       & 0.945  & 0.955  & 0.948  & 0.948  & 0.938  & 0.980  & \textbf{0.988} \\
			\textit{image3} &       & 0.2   & 27.98  & 17.98  & 27.92  & 27.92  & 29.42  & 32.65  & \textbf{33.74} &       & 0.921  & 0.428  & 0.932  & 0.932  & 0.941  & 0.978  & \textbf{0.984} \\
			& 0.6   & 0.1   & 28.88  & 29.96  & 28.45  & 28.46  & 31.02  & 33.60  & \textbf{36.08} &       & 0.953  & 0.959  & 0.952  & 0.952  & 0.961  & 0.983  & \textbf{0.991} \\
			&       & 0.2   & 28.31  & 29.34  & 28.05  & 28.05  & 29.67  & 33.20  & \textbf{34.70} &       & 0.933  & 0.944  & 0.941  & 0.941  & 0.948  & 0.980  & \textbf{0.986} \\
			& 0.8   & 0.1   & 29.52  & 30.28  & 28.81  & 28.81  & 33.29  & 34.79  & \textbf{38.44} &       & 0.964  & 0.965  & 0.957  & 0.957  & 0.970  & 0.986  & \textbf{0.995} \\
			&       & 0.2   & 29.08  & 29.80  & 28.44  & 28.44  & 31.27  & 34.19  & \textbf{36.63} &       & 0.951  & 0.957  & 0.950  & 0.950  & 0.961  & 0.984  & \textbf{0.992} \\
			\cline{1-10}\cline{12-18}    \textit{Kodak} & 0.5   & 0.1   & 29.75  & 29.91  & 28.22  & 28.23  & 30.16  & 31.15  & \textbf{32.45} &       & 0.950  & 0.950  & 0.942  & 0.942  & 0.952  & 0.964  & \textbf{0.977} \\
			\textit{image12} &       & 0.2   & 28.72  & 23.57  & 27.52  & 27.53  & 29.09  & 31.10  & \textbf{32.15} &       & 0.933  & 0.771  & 0.931  & 0.931  & 0.936  & 0.960  & \textbf{0.972} \\
			& 0.6   & 0.1   & 30.01  & 29.93  & 28.20  & 28.20  & 30.49  & 31.67  & \textbf{33.15} &       & 0.956  & 0.953  & 0.944  & 0.944  & 0.958  & 0.967  & \textbf{0.981} \\
			&       & 0.2   & 29.27  & 29.36  & 27.76  & 27.77  & 29.52  & 31.61  & \textbf{32.50} &       & 0.944  & 0.942  & 0.936  & 0.936  & 0.945  & 0.964  & \textbf{0.975} \\
			& 0.8   & 0.1   & 31.00  & 30.43  & 28.64  & 28.64  & 31.75  & 32.68  & \textbf{35.08} &       & 0.968  & 0.959  & 0.950  & 0.950  & 0.973  & 0.972  & \textbf{0.989} \\
			&       & 0.2   & 30.11  & 29.80  & 28.05  & 28.05  & 30.58  & 32.37  & \textbf{33.95} &       & 0.957  & 0.952  & 0.943  & 0.943  & 0.959  & 0.969  & \textbf{0.983} \\
			\hline
\end{tabular}}}%
\end{table}%
All numerical results are listed in Table \ref{tab:result_color_images} and the best PSNR and SSIM values are shown in bold.
We can see that NRQMC-NSS2D performs better than the other methods in all cases. Comparing the methods without NSS prior, the proposed NRQMC outperforms the others. Considering the NSS prior, PQMC-NSS2D and NRQMC-NSS2D can achieve better results than PQMC and NRQMC, and the NRQMC-NSS2D shows more improvement.
Figure \ref{fig:result_rgb1} shows the recovered results by different methods with $SR=0.2$ and $\gamma = 0.2$ for \textit{House}, \textit{Peppers} and \textit{Lena}.
Figure \ref{fig:result_rgb2} shows the recovered results by different methods with $SR=0.5$ and $\gamma = 0.1$ for \textit{Kodak image2}, \textit{Kodak image3} and \textit{Kodak image 12}.
It is seen that all methods can effectively restore the missing pixels of the images and remove most of the noise. We can observe that NRQMC produces clearer images and retains more details by comparing the results obtained from TNN, LRQA-G, QNN, PQMC, and NRQMC, such as the window in \textit{House}, the decoration on the hat in \textit{Lena}, and the text on the hat in \textit{Kodak image 3}.
However, it also includes some noises.
Adding PQMC-NSS2D and NRQMC-NSS2D to the comparison shows that the results from these two approaches contain less noise and are closer to the original.
The NRQMC-NSS2D preserves more details while the results obtained by PQMC-NSS2D are over-smooth, such as the window in \textit{House}, the pepper in \textit{Peppers} and the wood in \textit{Kodak image 2}.
It also can be seen that the results from NRQMC in Figure \ref{fig:result_rgb2} are visually better than those from PQMC-NSS2D.

\subsection{Color Videos}
In this subsection, we select four color videos to demonstrate the performance of the proposed methods, i.e., \textit{Akiyo}, \textit{Foreman}, \textit{Salesman} and \textit{Hall}\footnote{\href{http://trace.eas.asu.edu/yuv/index.html}{http://trace.eas.asu.edu/yuv/index.html}}. All videos are of $144 \times 176 \times 40$. We use the mean of PSNR (MPSNR) and the mean of SSIM (MSSIM) to evaluate the performance of the proposed model and algorithm.

We compare the proposed NRQMC-NSS3D for color video completion with LRQTC \cite{miao2020lowrank},  HaLRTC \cite{liu2013tensor}, SiLRTC \cite{liu2013tensor}, TMac-inc \cite{xu2015parallel}, TMac-dec \cite{xu2015parallel}, and PQMC-NSS3D. The PQMC-NSS3D uses PQMC to solve the subproblems in Step 3 of Algorithm \ref{alg4.1}. The $SR$ is set to be 0.15, 0.2, 0.3, 0.4, and 0.5.
The numerical results are listed in Table \ref{tab:result_color_videos_tc} and the best MPSNR and MSSIM values are shown in bold.
The NRQMC-NSS3D outperforms HaLRTC, SiLRTC, TMac-inc, TMac-dec, LRQTC, and PQMC-NSS3D in terms of MPSNR and MSSIM values.
Figure \ref{fig:result_cv_tc} shows the recovered results by different methods on the 25th frame of \textit{Akiyo}, \textit{Foreman}, \textit{Salesman}, and \textit{Hall} with $SR=0.3$.
All methods can effectively restore the missing pixels of the images. 
We can see that the proposed NRQMC-NSS3D outperforms HaLRTC, SiLRTC, TMac-inc, TMac-dec, and PQMC-NSS3D.

For robust completion, we compare the proposed NRQMC with PQMC under three different strategies: 1) recovering each frame of the video without NSS prior; 2) recovering each frame of the video using Algorithm \ref{alg4.1}; 3) recovering the color video using Algorithm \ref{alg4.1}.
Therefore, six methods are tested, i.e., PQMC, NRQMC, PQMC-NSS2D, NRQMC-NSS2D, PQMC-NSS3D, and NRQMC-NSS3D.
The result by NRQMC is used as the initial estimate of each NSS-based method.
The $SR$ is set to be 0.5, 0.6, and 0.8 and the $\gamma$ is set to be 0.1 and 0.2.
The numerical results are listed in Table \ref{tab:result_color_videos} and the best MPSNR and MSSIM values are shown in bold.
We can see that NRQMC-NSS3D achieves the highest MPSNR and MSSIM values in most cases.
The results obtained by NSS3D-based methods are better than those by NSS2D-based methods.
Under the same processing strategy, NRQMC performs better than PQMC in most cases.
Figure \ref{fig:result_color_videos} shows the recovered results by different methods on the 25th frame of \textit{Akiyo}, \textit{Foreman}, \textit{Salesman}, and \textit{Hall} with $SR=0.6$ and $\gamma = 0.1$.
It is seen that all methods can effectively restore the missing pixels of the images and remove most of the noise and NRQMC-NSS3D produces the best result. We can see that the methods using NRQMC can preserve more details while the methods using PQMC result in over-smoothness.

\begin{table}[H]
\centering
\caption{The MPSNR and MSSIM values of all methods for color videos.}\label{tab:result_color_videos_tc}
\scalebox{.52}{\setlength{\tabcolsep}{1.3mm}{

}
\caption{Recovery results by different methods on \textit{Akiyo}, \textit{Foreman}, \textit{Salesman}, and \textit{Hall} (from top to bottle) under $SR=0.3$. (a) Observed. (b) HaLRTC. (c) SiLRTC. (d) TMac-inc. (e) TMac-dec. (f) LRQTC. (g) PQMC-NSS3D. (h) NRQMC-NSS3D. (i) Original. }\label{fig:result_cv_tc}%
\end{figure}%
\begin{table}[H]
\centering
\caption{The MPSNR and MSSIM values of all methods for color videos.}\label{tab:result_color_videos}
\scalebox{.52}{\setlength{\tabcolsep}{1.3mm}{%
}
\caption{Recovery results by different methods on color videos under $SR=0.6$ and $\gamma=0.1$. (a) Observed. (b) PQMC. (c) NRQMC. (d) PQMC-NSS2D. (e) NRQMC-NSS2D. (f) PQMC-NSS3D. (g) NRQMC-NSS3D. (h) Original.}
\label{fig:result_color_videos}%
\end{figure}%

\section{Concluding remarks}
\label{sec:con}
In this paper, we present a new nonconvex approach for the robust QMC problem. We not only make use of
the nonconvex MCP approximation rank function under the QSVD framework for capturing the low-rank information but also introduce a suitable nonconvex sparsity measuration for the sparse constraint term rather than by the  $L_1$-norm commonly used in literature. The ADMM algorithm is established for solving
the proposed nonconvex robust QMC model, and its convergence analysis is given. In addition, the NSS prior is applied to search silimar patches, gathered into a low-rank quaternion matrix, of a color
image/video, and a new reconstruction is computed by the new nonconvex robust QMC algorithm. Numerical experiments on color image and video inpainting show the superiority of the proposed model and method, resulting in better PSNR and SSIM values.

Note that the proposed  model and algorithm still has several limitations. There is no advantage in the running time because it needs to compute QSVD at each iteration, which is time consuming for large-scale data.  In the future, we aim to   design better and faster QSVD approach to improve the efficiency of the model and algorithm. Besides, since we
mainly want to show the advantages of nonconvex surrogates under the quaternion framework, thus we do not assign weights on different singular
values and different sparse entries, which may cause underperformance in some special cases,
and it will be left for our future work. Besides, it is worth analyzing the improvement of the NRQMC compared with the QNN. Finally, in the future, we are interested  in studying the exact recovery condition and sampling size requirement of the NRQMC for successful recovery. 


\section*{Declaration of competing interest}
The authors declare that they have no known competing financial interests or personal
relationships that could have appeared to influence the work reported in this paper.

\bibliographystyle{elsarticle-num}
\bibliography{NRQMC}

\section*{Appendix. Proof of Main Results}

\subsection*{A1. Proof of Theorem \ref{thm3.2}}
From (\ref{3.14}), we see that  $g(\mathbf{x})$ is not differentiable at $\mathbf{x}=0$ and is continuous differentiable at $\mathbf{x}\neq0$. 
By direct calculations,  the first  order   gradient  of $g(\mathbf{x})$ is given by
\begin{equation}\label{3.15}
\nabla g(\mathbf{x})=\mathbf{x}-\mathbf{y}+\nu p|\mathbf{x}|^{p-1}\mathrm{signQ}(\mathbf{x}),
\end{equation}
where $\mathbf{x}\neq0$.
For $0\neq\mathbf{x}\in\mathbb{Q}$ and $\bm{\delta}\mathbf{x}\in\mathbb{Q}$, we have
\begin{small}
\begin{eqnarray}
&&\nabla g(\mathbf{x}+\bm{\delta}\mathbf{x}) -\nabla g(\mathbf{x})\nonumber \\
&=&[\mathbf{x}+\bm{\delta}\mathbf{x}-\mathbf{y}+\nu p|\mathbf{x}+\bm{\delta}\mathbf{x}|^{p-1}\mathrm{signQ}(\mathbf{x}+\bm{\delta}\mathbf{x})]-[\mathbf{x}-\mathbf{y}+\nu p|\mathbf{x}|^{p-1}\mathrm{signQ}(\mathbf{x})]\nonumber \\
&=&\bm{\delta}\mathbf{x}+\nu p|\mathbf{x}+\bm{\delta}\mathbf{x}|^{p-2}(\mathbf{x}+\bm{\delta}\mathbf{x})-\nu p|\mathbf{x}|^{p-2} \mathbf{x} \nonumber \\
&=&\bm{\delta}\mathbf{x}+\nu p(|\mathbf{x}|^2+\langle\mathbf{x},\bm{\delta}\mathbf{x}\rangle+|\bm{\delta}\mathbf{x}|^2)^{\frac{p-2}{2}}(\mathbf{x}+\bm{\delta}\mathbf{x})-\nu p|\mathbf{x}|^{p-2} \mathbf{x} \nonumber \\
&=&\bm{\delta}\mathbf{x}+\nu p|\mathbf{x}|^{p-2}\Big(1+2\frac{\langle\mathbf{x},\bm{\delta}\mathbf{x}\rangle}{|\mathbf{x}|^2}
+\frac{|\bm{\delta}\mathbf{x}|^2}{|\mathbf{x}|^2}\Big)^{\frac{p-2}{2}}(\mathbf{x}+\bm{\delta}\mathbf{x})-\nu p|\mathbf{x}|^{p-2} \mathbf{x} \nonumber \\
&=&\bm{\delta}\mathbf{x}+\nu p|\mathbf{x}|^{p-2}\Big(1+(p-2)\frac{\langle\mathbf{x},\bm{\delta}\mathbf{x}\rangle}{|\mathbf{x}|^2}
+o(|\bm{\delta}\mathbf{x}|)\Big) (\mathbf{x}+\bm{\delta}\mathbf{x})-\nu p|\mathbf{x}|^{p-2} \mathbf{x} \nonumber \\
&=&(1+\nu p|\mathbf{x}|^{p-2})\bm{\delta}\mathbf{x}
+\nu p(p-2)|\mathbf{x}|^{p-2} \langle\mathrm{signQ}(\mathbf{x}),\bm{\delta}\mathbf{x}\rangle \mathrm{signQ}(\mathbf{x})+o(|\bm{\delta}\mathbf{x}|).
\end{eqnarray}
\end{small}
This means that
\begin{equation*}
\nabla^2g(\mathbf{x}) \bm{\delta}\mathbf{x} =(1+\nu p|\mathbf{x}|^{p-2})\bm{\delta}\mathbf{x}
+\nu p(p-2)|\mathbf{x}|^{p-2} \langle\mathrm{signQ}(\mathbf{x}),\bm{\delta}\mathbf{x}\rangle \mathrm{signQ}(\mathbf{x}).
\end{equation*}
Then
\begin{small}
\begin{eqnarray*}
&&\langle \nabla^2g(\mathbf{x}) \bm{\delta}\mathbf{x},\bm{\delta}\mathbf{x}\rangle
=(1+\nu p|\mathbf{x}|^{p-2})|\bm{\delta}\mathbf{x}|^2
+\nu p(p-2)|\mathbf{x}|^{p-2} \langle\mathrm{signQ}(\mathbf{x}),\bm{\delta}\mathbf{x}\rangle^2>0
\end{eqnarray*}
\end{small}
for arbitrary $\bm{\delta}\mathbf{x}\neq 0$.
Thus $\nabla^2 g (\mathbf{x})$ is positive definite  when $\mathbf{x}\neq 0$. Let $\nabla g(\mathbf{x}_\clubsuit)=0$, i.e.,
$
\mathbf{x}_\clubsuit-\mathbf{y}+\nu p|\mathbf{x}_\clubsuit|^{p-1}\mathrm{signQ}(\mathbf{x}_\clubsuit)=0$.
It follows from Proposition \ref{prop2.2}  that $\mathbf{x}_\clubsuit$ is an extremum point of $g$. Correspondingly, $\max\{g(\mathbf{x}_\clubsuit),g(0)\}$ is the minimum value of $g$.  There may exist  a specific
$\mathbf{y}$, where $g(\mathbf{x}_\clubsuit)$ is exactly $g(0)$. Thus, to generalize soft thresholding, we should solve the following nonlinear equation system to determine a correct thresholding value
$\tau_p^{\mathrm{QGST}}(\nu)$ ($|\mathbf{y}|$) and its corresponding $\mathbf{x}_\clubsuit$:
\begin{subequations}
\begin{align}
	\frac{1}{2}|\mathbf{x}_\clubsuit-\mathbf{y}|^2+\nu|\mathbf{x}_\clubsuit|^p
	=\frac{1}{2}|\mathbf{y}|^2\label{3.17a}\\
	\mathbf{x}_\clubsuit-\mathbf{y}+\nu p|\mathbf{x}_\clubsuit|^{p-1}\mathrm{signQ}(\mathbf{x}_\clubsuit)=0.\label{3.17b}
\end{align}
\end{subequations}
By (\ref{3.17b}), it follows that
\begin{equation}\label{3.18}
|\mathbf{y}|=|\mathbf{x}_\clubsuit|+\nu p|\mathbf{x}_\clubsuit|^{p-1}~~
\text{and}~~
\mathbf{y}=\mathbf{x}_\clubsuit+\nu p|\mathbf{x}_\clubsuit|^{p-1}\mathrm{signQ}(\mathbf{x}_\clubsuit).
\end{equation}
Substituting the above equalities into (\ref{3.17a}) yields
\begin{equation*}
\frac{1}{2}|\nu p|\mathbf{x}_\clubsuit|^{p-1}\mathrm{signQ}(\mathbf{x}_\clubsuit)|^2
+\nu|\mathbf{x}_\clubsuit|^p=\frac{1}{2}(|\mathbf{x}_\clubsuit|+\nu p|\mathbf{x}_\clubsuit|^{p-1})^2,
\end{equation*}
and then
\begin{equation}\label{3.19}
|\mathbf{x}_\clubsuit|^p(2\nu(1-p)-|\mathbf{x}_\clubsuit|^{2-p})=0.
\end{equation}
This means that $|\mathbf{x}_\clubsuit|=(2\nu (1-p))^{\frac{1}{2-p}}$. This, together with  (\ref{3.18}), yields
\begin{equation}\label{3.20}
\tau_p^{\mathrm{QGST}}(\nu)=|\mathbf{y}|=|\mathbf{x}_\clubsuit|+\nu p|\mathbf{x}_\clubsuit|^{p-1}=(2\nu (1-p))^{\frac{1}{2-p}}+\nu p(2\nu (1-p))^{\frac{p-1}{2-p}}.
\end{equation}
Now the last thing is to determine $\mathbf{x}_\clubsuit$. By (\ref{3.17b}), one has
\begin{eqnarray}
\mathbf{x}_\clubsuit &=& \frac{|\mathbf{x}_\clubsuit|}{|\mathbf{x}_\clubsuit|+\nu p|\mathbf{x}_\clubsuit|^{p-1}}\mathbf{y}=\frac{|\mathbf{x}_\clubsuit|}{|\mathbf{y}|}\mathbf{y}
=\mathrm{signQ}(\mathbf{y})|\mathbf{x}_\clubsuit|.\label{3.21}
\end{eqnarray}
Fortunately, we see from the first equality of (\ref{3.18}) that $|\mathbf{x}_\clubsuit|$ can be determined by solving  $|\mathbf{y}|=|\mathbf{x}_\clubsuit|+\nu p|\mathbf{x}_\clubsuit|^{p-1}$ using fixed point iteration. Once $|\mathbf{x}_\clubsuit|$ is obtained, by (\ref{3.21}), $\mathbf{x}_\clubsuit$ is achieved. The proof is completed.

\subsection*{A2. Proof of Theorem \ref{thm3.1}}
Since $\mathbf{Y}_k=\mathbf{U}_k\mathbf{\Sigma}_k\mathbf{V}_k^*$, it follows that $\mathbf{\Sigma}_k=\mathbf{U}_k^*\mathbf{Y}_k\mathbf{V}_k$. Using the definition of MCP function in (\ref{3.2}), we have
\begin{eqnarray}
\min\frac{1}{\mu_k}\|\mathbf{L}\|_{\text{MCP}}+
\frac{1}{2}\|\mathbf{L}-\mathbf{Y}_k\|_F^2
\sum\limits_{i=1}^{n_{(2)}}\frac{1}{\mu_k}\Phi_{c,\eta}(\sigma_i(\mathbf{L}))
+\frac{1}{2}\|\mathbf{U}_k^*\mathbf{L}\mathbf{V}_k-\mathbf{\Sigma}_k\|_F^2,\label{3.11}
\end{eqnarray}
where the second equality comes from the unitary invariance of quaternion Frobenius norm.
Set $\mathbf{A}=\mathbf{U}_k^*\mathbf{L}\mathbf{V}_k$. Then $\mathbf{A}$ and $\mathbf{L}$ have the same singular values. Assume that $\mathbf{\Sigma}_\mathbf{L}$ and $\mathbf{\Sigma}_\mathbf{A}$ are the singular value matrices of $\mathbf{L}$ and $\mathbf{A}$, respectively. Thus $\mathbf{\Sigma}_\mathbf{L}=\mathbf{\Sigma}_\mathbf{A}$  and Eq.(\ref{3.11}) is equal to
{\footnotesize\begin{eqnarray}
	\!\!\!\!\!\!\!\!&&\sum\limits_{i=1}^{n_{(2)}}\frac{1}{\mu_k}\Phi_{c,\eta}(\sigma_i(\mathbf{L}))
	+\frac{1}{2}\|\mathbf{A}-\mathbf{\Sigma}_k\|_F^2
	\geq \sum\limits_{i=1}^{n_{(2)}}\frac{1}{\mu_k}\Phi_{c,\eta}(\sigma_i(\mathbf{L}))
	+\frac{1}{2}\|\mathbf{\Sigma}_\mathbf{A}-\mathbf{\Sigma}_k\|_F^2\nonumber\\
	\!\!\!\!\!\!\!\!&=&\!\!\!\! \sum\limits_{i=1}^{n_{(2)}}\frac{1}{\mu_k}\Phi_{c,\eta}(\sigma_i(\mathbf{L}))
	\!+\!\frac{1}{2}\|\mathbf{\Sigma}_\mathbf{L}\!-\!\mathbf{\Sigma}_k\|_F^2
	\!	=\!\sum\limits_{i=1}^{n_{(2)}}\Big(\frac{1}{\mu_k}\Phi_{c,\eta}(\sigma_i(\mathbf{L}))
	\!+\!\frac{1}{2}\|\sigma_i(\mathbf{L})\!-\!(\mathbf{\Sigma}_k)_{i,i}\|^2\Big) .\label{3.12}
\end{eqnarray}}
where the first inequality follows from the von Neumann trace theorem (\cite[Proposition 1]{jia2019robust}). Hence, solving the original minimization problem in Eq.(\ref{3.10}) is transformed to solving Eq.(\ref{3.12}) and the optimal solution $\mathbf{L}_\clubsuit$ is $\mathbf{L}_\clubsuit=\mathbf{U}_k\mathrm{Diag}\Big(
\mathrm{Prox}_{\frac{1}{\mu_k}\Phi_{c,\eta}}\big((\mathbf{\Sigma}_k)_{1,1}\big),\cdots,
\mathrm{Prox}_{\frac{1}{\mu_k}\Phi_{c,\eta}}\big((\mathbf{\Sigma}_k)_{n_{(2)},n_{(2)}}\big) \mathbf{V}_k^*$. This obtains the desired result.

\subsection*{A3. Proof of Theorem \ref{thm4.1}}

\begin{lem}\label{lem4.2}
The sequence $\{\mathbf{M}_k\}$ generated by {\rm Algorithm \ref{alg3.3}} is bounded.
\end{lem}
{\bf Proof.}~By the definition of MCP quaternion rank approximation, we have
{\footnotesize	\begin{equation}\label{4.2}
\nabla \|\mathbf{L}\|_{\mathrm{MCP}}\!=\!
\mathbf{U}\mathrm{Diag} \big(
\partial\Phi_{c,\eta}(\sigma_1(\mathbf{L}))/\partial\sigma_1(\mathbf{L}),\cdots,
\partial\Phi_{c,\eta}(\sigma_{n_{(2)}}(\mathbf{L}))/\partial\sigma_{n_{(2)}}(\mathbf{L})\big)
\mathbf{V}^*.
\end{equation}}
It is easily seen that
$\partial\Phi_{c,\eta}(\sigma_i(\mathbf{L}))/\partial\sigma_i(\mathbf{L})\leq c,~
\forall 1\leq i\leq n_{(2)}$.
By (\ref{4.2}),  it follows that $\|\nabla \|\mathbf{L}\|_{\mathrm{MCP}}\|_F^2\leq n_{(2)}c^2$. Therefore, $\nabla \|\mathbf{L}\|_{\mathrm{MCP}}$ is bounded.

From the first order optimality condition of (\ref{3.10})  in $\mathbf{L}_{k+1}$, we have
\begin{equation*}
\frac{1}{\mu_k}\nabla \|\mathbf{L}\|_{\mathrm{MCP}}|_{\mathbf{L}=\mathbf{L}_{k+1}}+\mathbf{L}_{k+1}-\mathbf{Y}_k=0,
\end{equation*}
where $\mathbf{Y}_k=\mathbf{X}-\mathbf{S}_{k+1}-\frac{\mathbf{M}_k}{\mu_k}$. Correspondingly, it holds
\begin{equation*}
\nabla \|\mathbf{L}\|_{\mathrm{MCP}}|_{\mathbf{L}=\mathbf{L}_{k+1}}
+[\mathbf{M}_k+\mu_k(\mathbf{L}_{k+1}
+\mathbf{S}_{k+1}-\mathbf{X})]=0.
\end{equation*}
From line 5 of Algorithm \ref{alg3.3}, it follows that
$\nabla \|\mathbf{L}\|_{\mathrm{MCP}}|_{\mathbf{L}=\mathbf{L}_{k+1}}
+\mathbf{M}_{k+1}=0$.
Thus  $\{\mathbf{M}_{k}\}$ appears to be bounded. The proof is completed. \hfill$\Box$

\begin{lem}\label{lem4.3}
The sequences $\{\mathbf{L}_k\}$ and $\{\mathbf{S}_k\}$ generated by {\rm Algorithm \ref{alg3.3}} are bounded.
\end{lem}
{\bf Proof.}~By the definition of the augmented Lagrange
function in (\ref{3.8}) and line 5 of Algorithm \ref{alg3.3}, it follows that
{\scriptsize\begin{eqnarray}
	&&\mathscr{L}_{\mu_k}(\mathbf{L}_k,\mathbf{S}_k,\mathbf{M}_k)\nonumber\\
	&=&\mathscr{L}_{\mu_{k-1}}(\mathbf{L}_k,\mathbf{S}_k,\mathbf{M}_{k-1})
	+\frac{\mu_k}{2}\Big\|\frac{\mathbf{M}_k-\mathbf{M}_{k-1}}{\mu_{k-1}}+\frac{\mathbf{M}_k}{\mu_k}\Big\|_F^2
	-\frac{\mu_{k-1}}{2}\Big\|\frac{\mathbf{M}_k}{\mu_{k-1}}
	\Big \|_F^2-\frac{\|\mathbf{M}_k\|_F^2}{2\mu_k}
	+\frac{\|\mathbf{M}_{k-1}\|_F^2}{2\mu_{k-1}}\nonumber\\
	&=&\mathscr{L}_{\mu_{k-1}}(\mathbf{L}_k,\mathbf{S}_k,\mathbf{M}_{k-1})
	+\frac{\mu_k}{2\mu_{k-1}^2}\|\mathbf{M}_k-\mathbf{M}_{k-1}\|^2_F+\frac{\|\mathbf{M}_k\|_F^2}{2\mu_k}
	+\frac{\|\mathbf{M}_k\|_F^2}{\mu_{k-1}}
	-\frac{1}{\mu_{k-1}}\langle\mathbf{M}_{k-1},\mathbf{M}_k\rangle\nonumber\\
	&&
	-\frac{\|\mathbf{M}_k\|_F^2}{2\mu_{k-1}}
	-\frac{1}{2\mu_k}\|\mathbf{M}_k\|_F^2
	+\frac{1}{2\mu_{k-1}}\|\mathbf{M}_{k-1}\|_F^2\nonumber\\
	&=&\mathscr{L}_{\mu_{k-1}}(\mathbf{L}_k,\mathbf{S}_k,\mathbf{M}_{k-1})
	+\frac{\mu_k+\mu_{k-1}}{2\mu_{k-1}^2}\|\mathbf{M}_k-\mathbf{M}_{k-1}\|^2_F.\label{4.3}
\end{eqnarray}}
Since
$\mathbf{S}_{k+1}=\arg\min_{\mathbf{S}} \mathscr{L}_{\mu_{k}}(\mathbf{L}_k,\mathbf{S},\mathbf{M}_{k})~~\text{and}~~
\mathbf{L}_{k+1}=\arg\min_{\mathbf{L}}\mathscr{L}_{\mu_{k}}(\mathbf{L},\mathbf{S}_{k+1},\mathbf{M}_{k})$,
we have
{\scriptsize\begin{eqnarray}
&&	\mathscr{L}_{\mu_{k}}(\mathbf{L}_{k+1},\mathbf{S}_{k+1},\mathbf{M}_{k})\leq   \mathscr{L}_{\mu_{k}}(\mathbf{L}_{k},\mathbf{S}_{k+1},\mathbf{M}_{k})
\leq\mathscr{L}_{\mu_{k}}(\mathbf{L}_{k},\mathbf{S}_{k},\mathbf{M}_{k})\nonumber\\
&\leq&\mathscr{L}_{\mu_{k-1}}(\mathbf{L}_k,\mathbf{S}_k,\mathbf{M}_{k-1})
+\frac{\mu_k+\mu_{k-1}}{2\mu_{k-1}^2}\|\mathbf{M}_k-\mathbf{M}_{k-1}\|^2_F
\leq  \mathscr{L}_{\mu_{0}}(\mathbf{L}_1,\mathbf{S}_1,\mathbf{M}_{0})
+\sum\limits_{i=1}^k\frac{\mu_i+\mu_{i-1}}{2\mu_{i-1}^2}\|\mathbf{M}_i-\mathbf{M}_{i-1}\|^2_F \nonumber\\
&\leq&  \mathscr{L}_{\mu_{0}}(\mathbf{L}_0,\mathbf{S}_1,\mathbf{M}_{0})
+\sum\limits_{i=1}^k\frac{\mu_i+\mu_{i-1}}{2\mu_{i-1}^2}\|\mathbf{M}_i-\mathbf{M}_{i-1}\|^2_F
\leq \frac{\mu_0}{2}\|\mathbf{X}\|_F^2
+\big(\max\limits_{i}\|\mathbf{M}_i-\mathbf{M}_{i-1}\|^2_F\big)
\sum\limits_{i=1}^k\frac{\mu_i+\mu_{i-1}}{2\mu_{i-1}^2}.\nonumber
\end{eqnarray}}
According to  Lemma \ref{lem4.2}, it follows that $\{\mathbf{M}_k\}$ is bounded. Then $\max\limits_{i}\|\mathbf{M}_i-\mathbf{M}_{i-1}\|^2_F$ is also bounded. From   Algorithm \ref{alg3.3}, we have $\mu_i=\rho\mu_{i-1}=\rho^i\mu_0$, $\rho=1.1$ and $\mu_0=10^{-4}$. This gives
\begin{equation*}
\sum\limits_{i=1}^\infty\frac{\mu_i+\mu_{i-1}}{2\mu_{i-1}^2}=
\frac{\rho+1}{2\mu_0}\sum\limits_{i=1}^\infty\frac{1}{\rho^{i-1}}
=\frac{\rho+1}{2\mu_0}.
\end{equation*}
Hence $\sum\limits_{i=1}^k\frac{\mu_i+\mu_{i-1}}{2\mu_{i-1}^2}$ is bounded, and then $\mathscr{L}_{\mu_{k}}(\mathbf{L}_{k+1},\mathbf{S}_{k+1},\mathbf{M}_{k})$ has upper bound.

Using the relation (\ref{3.8}) again, we have
\begin{eqnarray}
&&\mathscr{L}_{\mu_{k}}(\mathbf{L}_{k+1},\mathbf{S}_{k+1},\mathbf{M}_{k})
+\frac{\|\mathbf{M}_k\|_F^2}{ 2\mu_k}\nonumber\\
&=&\|\mathbf{L}_{k+1}\|_{\text{MCP}}+
\lambda \|\mathcal{P}_\Omega(\mathbf{S}_{k+1})\|_p^p+ \frac{\mu_k}{2}\Big\|\mathbf{L}_{k+1}+\mathbf{S}_{k+1}-\mathbf{X}+\frac{\mathbf{M}_k}{\mu_k}\Big\|_F^2.
\label{4.6}
\end{eqnarray}
Since every term on the right-hand side of the equation (\ref{4.6}) is nonnegative, it follows that the sequences $\{\mathbf{L}_{k+1}\}$ and $\{\mathbf{S}_{k+1}\}$ are bounded. The proof is completed.\hfill$\Box$

\begin{lem}\label{lem4.4}
For the quaternion $L_p$-regularized unconstrained optimization problem
\begin{equation}\label{4.7}
	\min\limits_{\mathbf{X}\in \mathbb{Q}^{n_1\times n_2}}\{F(\mathbf{X}):=f(\mathbf{X})+\lambda\|\mathbf{X}\|_p^p\},
\end{equation}
where $\lambda>0$, $p\in(0, 1)$, and $f$ is continuous differentiable in $\mathbb{Q}^{n_1\times n_2}$.
Let $\mathbf{X}_\clubsuit$ be a local minimizer of {\rm(\ref{4.7})}. Then $\mathbf{X}_\clubsuit$ is a first order
stationary point, i.e., $\langle\mathbf{X}_\clubsuit,\nabla f(\mathbf{X}_\clubsuit)\rangle + \lambda p \|\mathbf{X}_\clubsuit\|_p^p = 0$ holds at $\mathbf{X}_\clubsuit$.
\end{lem}
{\bf Proof.}~Let $\mathbf{X}_\clubsuit=(\mathbf{x}_{ij})\in \mathbb{Q}^{n_1\times n_2}$.
By the definition of quaternion $L_p$-norm and some direct calculations,  we have
\begin{equation*}
\partial f(\mathbf{X})/\partial\mathbf{x}_{ij}+\lambda p|\mathbf{x}_{ij}|^{p-1}\mathrm{signQ}(\mathbf{x}_{ij})=0,i=1,\cdots,n_1,j=1,\cdots,n_2.
\end{equation*}
Multiplying $\mathbf{x}_{ij}^*$ from the left of the above equality yields
\begin{equation*}
\mathbf{x}_{ij}^*\partial f(\mathbf{X})/\partial\mathbf{x}_{ij}+\lambda p|\mathbf{x}_{ij}|^{p}=0,i=1,\cdots,n_1,j=1,\cdots,n_2.
\end{equation*}
This means that
$	\sum\limits_{i=1}^{n_1}\sum\limits_{j=1}^{n_2}\Big(\mathbf{x}_{ij}^*\partial f(\mathbf{X})/\partial\mathbf{x}_{ij}+\lambda p|\mathbf{x}_{ij}|^{p}\Big)=0$,
and then $\langle\mathbf{X}_\clubsuit$, $\nabla f(\mathbf{X}_\clubsuit)\rangle + \lambda p \|\mathbf{X}_\clubsuit\|_p^p = 0$. The proof is completed. \hfill$\Box$

\medskip
Based on Lemmas \ref{lem4.2}-\ref{lem4.4}, we give the proof of Theorem \ref{thm4.1}.

\medskip

By Lemmas \ref{lem4.2} and \ref{lem4.3},  the sequence $\mathbf{W}_k = \{\mathbf{L}_k, \mathbf{S}_k,\mathbf{M}_k\}$  generated by Algorithm \ref{alg3.3} is bounded.
From the Bolzano-Weierstrass theorem, the sequence  $\{\mathbf{W}_k\}$  has at
least one convergent subsequence, thus there exists at least one accumulation point $\mathbf{W}_\clubsuit = \{\mathbf{L}_\clubsuit, \mathbf{S}_\clubsuit,\mathbf{M}_\clubsuit\}$ of $\mathbf{W}_k$ for the sequence $\{\mathbf{W}_k\}$. Without loss of generality, we assume that $\mathbf{W}_k$ converges
to $\mathbf{W}_\clubsuit$.

According to line 5 of Algorithm \ref{alg3.1}, it follows that
\begin{equation*}
\lim\limits_{k\rightarrow\infty} (\mathbf{L}_{k+1}+\mathbf{S}_{k+1}-\mathbf{X})=\lim\limits_{k\rightarrow\infty} \frac{\mathbf{M}_{k+1}-\mathbf{M}_k}{\mu_k}=0.
\end{equation*}
This implies that $\mathbf{L}_\clubsuit+\mathbf{S}_\clubsuit=\mathbf{X}$.
Since
\begin{equation*}
\mathbf{L}_{k+1}=\arg\min\limits_{\mathbf{L}}\mathscr{L}_{\mu_{k}}(\mathbf{L},\mathbf{S}_{k+1},\mathbf{M}_{k})
=\arg\min\limits_{\mathbf{L}}\|\mathbf{L}\|_{\text{MCP}}+
\frac{\mu_k}{2}\Big\|\mathbf{L}+\mathbf{S}_{k+1}-\mathbf{X}+\frac{\mathbf{M}_k}{\mu_k}\Big\|_F^2,
\end{equation*}
we have
{\footnotesize\begin{eqnarray*}
	0=\nabla \|\mathbf{L}\|_{\text{MCP}}|_{\mathbf{L}=\mathbf{L}_{k+1}}+
	\mu_k\Big(\mathbf{L}_{k+1}+\mathbf{S}_{k+1}-\mathbf{X}+\frac{\mathbf{M}_k}{\mu_k}\Big)
	=\nabla \|\mathbf{L}\|_{\text{MCP}}|_{\mathbf{L}=\mathbf{L}_{k+1}}+\mathbf{M}_{k+1}.
\end{eqnarray*}}
This implies that	$\nabla \|\mathbf{L}\|_{\mathrm{MCP}}|_{\mathbf{L}=\mathbf{L}_\clubsuit}+\mathbf{M}_\clubsuit=0$.
Since
{\footnotesize \begin{equation*}
	\mathscr{P}_\Omega(\mathbf{S}_{k+1})=\arg\min\limits_{\mathbf{S}}
	\mathscr{L}_{\mu_{k}}(\mathbf{L}_{k},\mathbf{S},\mathbf{M}_{k})
	=\arg\min\limits_{\mathbf{S}}\lambda \|\mathcal{P}_\Omega(\mathbf{S})\|_p^p+ \frac{\mu_k}{2}\Big\|\mathcal{P}_\Omega\big(\mathbf{L}_{k}+\mathbf{S}-\mathbf{X}
	+\frac{\mathbf{M}_k}{\mu_k}\big)\Big\|_F^2,
\end{equation*}}
it follows from Lemma \ref{lem4.4} that
\begin{equation*}
\lambda \|\mathcal{P}_\Omega(\mathbf{S}_{k+1})\|_p^p+ \big\langle\mathcal{P}_\Omega(\mathbf{S}_{k+1}),\mu_k\mathcal{P}_\Omega\big(\mathbf{L}_{k}+\mathbf{S}_{k+1}-\mathbf{X}
+\frac{\mathbf{M}_k}{\mu_k}\big)\big\rangle=0.
\end{equation*}
This, together with line 5 of Algorithm \ref{alg3.3}, yields
\begin{equation}\label{4.11}
\lambda \|\mathcal{P}_\Omega(\mathbf{S}_{k+1})\|_p^p+ \langle\mathcal{P}_\Omega(\mathbf{S}_{k+1}),\mathcal{P}_\Omega(
\mathbf{M}_{k+1}-\mu_k(\mathbf{L}_{k+1}-\mathbf{L}_k))\rangle=0.
\end{equation}
Notice that $\mu_k <\mu_{\max}= 10^{10}$, we have $\mu_k (\mathbf{L}_{k+1}-\mathbf{L}_k)\rightarrow 0$ if $k\rightarrow\infty$. From
(\ref{4.11}), it follows that
$\lambda \|\mathcal{P}_\Omega(\mathbf{S}_\clubsuit)\|_p^p+ \langle\mathcal{P}_\Omega(\mathbf{S}_\clubsuit),\mathcal{P}_\Omega(
\mathbf{M}_\clubsuit)\rangle=0$.
Since
\begin{equation*}
\mathscr{P}_{\Omega^\perp}(\mathbf{S}_{k+1})=\arg\min\limits_{\mathbf{S}}
\mathscr{L}_{\mu_{k}}(\mathbf{L}_{k},\mathbf{S},\mathbf{M}_{k}) =\arg\min_{\mathbf{S}}\frac{\mu_k}{2}\Big\|\mathcal{P}_{\Omega^\perp}\big(\mathbf{L}_{k}+\mathbf{S}-\mathbf{X}
+\frac{\mathbf{M}_k}{\mu_k}\big)\Big\|_F^2,
\end{equation*}
we have
$
0= \mu_k\mathcal{P}_{\Omega^\perp}\Big(\mathbf{L}_{k}+\mathbf{S}_{k+1}-\mathbf{X}
+\frac{\mathbf{M}_k}{\mu_k}\Big)=\mathcal{P}_{\Omega^\perp}(
\mathbf{M}_{k+1}-\mu_k(\mathbf{L}_{k+1}-\mathbf{L}_k))$.
Letting $k\rightarrow\infty$, we obtain $\mathcal{P}_{\Omega^\perp}(
\mathbf{M}_\clubsuit)=0$,
Thus,   the desired result is obtained.

\end{document}